\documentclass[reqno]{amsart}

\usepackage{amsmath}
\usepackage{amssymb}
\usepackage{verbatim}
\usepackage{dsfont}
\usepackage{enumerate}
\usepackage[all]{xy}
\usepackage[mathscr]{eucal}
\usepackage{amsthm}
\usepackage{stmaryrd}
\usepackage{amscd}
\usepackage{amsfonts}
\usepackage{latexsym}
\usepackage{amscd}
\usepackage{graphicx}
\usepackage{tikz}
\usepackage{setspace}
\usepackage{txfonts}
\usepackage{wasysym}
\usepackage{todonotes}
\usepackage{hyperref}
\usetikzlibrary{decorations.pathmorphing}

\makeatletter

\newcommand{\arxiv}[1]{\href{http://arxiv.org/abs/#1}{\tt
    arXiv:\nolinkurl{#1}}}

\def\dot{\color{darkblue}{\color{white}\bullet}\!\!\!\circ}
\def\diamond{\color{white}\scriptstyle\spadesuit\color{darkblue}\hspace{-1.328mm}\scriptstyle\varspadesuit}

\def\bd#1{\text{\boldmath${#1}$}}
\newcommand{\jontodo}{\todo[inline,color=green!20]}

\newtheorem{theorem}{Theorem}[section]
\newtheorem{lemma}[theorem]{Lemma}

\theoremstyle{definition}  
\newtheorem{definition}[theorem]{Definition}

\newtheorem{remark}[theorem]{Remark}

\def\E{\up}
\def\F{\down}
\def\op{{\operatorname{op}}}

\@addtoreset{equation}{section}
\def\theequation{\arabic{section}.\arabic{equation}}

\def\coind{\operatorname{Coind}}
\def\ind{\operatorname{Ind}}
\def\res{\operatorname{Res}}

\def\Sym{\operatorname{Sym}}
\def\Pol{\operatorname{Pol}}
\def\H{\mathcal{H}eis}
\newcommand{\Seq}{\operatorname{Seq}}
\newcommand{\SEQ}{\operatorname{SEQ}}
\newcommand{\End}{\operatorname{End}}
\newcommand{\wt}{\operatorname{wt}}
\newcommand{\Aut}{\operatorname{End}} 
\newcommand{\CatHom}{\mathcal{H}om}
\newcommand{\unit}{\mathds{1}}
\def\ev{\operatorname{Ev}}
\def\up{{\color{darkblue}\uparrow}}
\def\down{{\color{darkblue}\downarrow}}
\newcommand{\Hom}{\operatorname{Hom}} 
\newcommand{\HOM}{\operatorname{HOM}} 
\newcommand{\Mod}{\operatorname{\!-mod}}
\newcommand{\proj}{\operatorname{\!-pmod}}
\newcommand{\Rep}{\operatorname{Rep}}
\newcommand{\id}{\text{Id}}
\newcommand{\gl}{\mathfrak{gl}}
\newcommand{\g}{\mathfrak{g}}
\newcommand{\Z}{\mathbb{Z}}
\newcommand{\N}{\mathbb{N}}
\newcommand{\C}{\mathbb{C}}
\newcommand{\Q}{\mathbb{Q}}
\newcommand{\K}{\mathbb{K}}
\newcommand{\eps}{\varepsilon}
\newcommand{\gr}{\operatorname{gr}}
\newcommand{\rev}{^{\operatorname{rev}}}
\newcommand{\rad}{{\operatorname{rad}}}
\def\la{\lambda}
\def\al{\alpha}
\def\be{\beta}
\def\Kar{\operatorname{Kar}}
\def\ga{\gamma}
\def\Ga{\Gamma}
\def\De{\Delta}

\def\dotM{\text{\boldmath$M$}}
\def\dotP{\text{\boldmath$P$}}
\def\dotp{\text{\boldmath$p$}}
\def\dotv{\text{\boldmath$v$}}
\def\dotk{\text{\boldmath$k$}}
\def\dote{\text{\boldmath$e$}}
\def\dotf{\text{\boldmath$f\!$}}

\newcommand{\Bim}{\mathcal{B}im}
\newcommand{\KB}{\mathcal{KB}}
\newcommand{\AOB}{\mathcal{AOB}}
\newcommand{\OB}{{\mathcal{OB}}}
\newcommand{\dotOB}{\dot{\mathcal O}{\mathcal B}}
\renewcommand{\Vec}{{\mathcal{V}ec}}
\renewcommand{\vec}{{\mathcal{V}ec}_{fd}}

\renewcommand{\k}{\Bbbk}

\definecolor{darkblue}{HTML}{111199}
\definecolor{darkgreen}{HTML}{336633}
\definecolor{darkred}{HTML}{993333}
\definecolor{darkpurple}{HTML}{995599}

\allowdisplaybreaks

\begin{document}

\title[Heisenberg categories]{On the definition of 
  Heisenberg category}

\author[J. Brundan]{Jonathan Brundan}
\address{Department of Mathematics,
University of Oregon, Eugene, OR 97403, USA}
\email{brundan@uoregon.edu}

\thanks{Research supported in part by NSF grant DMS-1700905.}

\begin{abstract}
We revisit the definition of 
the Heisenberg category of central charge $k \in \Z$.
For central charge $- 1$, this category was
introduced originally by Khovanov, but with
some additional cyclicity relations which we show here are
unnecessary.
For other negative central charges, the definition is due to Mackaay and
Savage, also with some redundant relations,
while central charge zero recovers the affine oriented Brauer
category of Brundan, Comes, Nash and Reynolds.
We also discuss cyclotomic quotients.
\end{abstract}

\maketitle  

\section{Introduction}

In \cite{K}, Khovanov introduced a graphical calculus for the
induction and restriction functors $\ind_{n}^{n+1}$ and $\res_{n-1}^{n}$
arising in the representation theory of the symmetric group $S_n$. This led
him to the definition of a monoidal category $\mathcal H$,
which he
called the {\em Heisenberg category}. This category is monoidally
generated by two objects $\up$ and $\down$ (corresponding to the induction and
restriction functors) with
morphisms defined in terms of equivalence classes of certain diagrams
modulo Reidemeister-type relations plus a small number of additional
relations. 
Khovanov's relations imply in particular that there is an {\em isomorphism} 
$$
\left[
\mathord{
\begin{tikzpicture}[baseline = -.4mm]
	\draw[<-,thick,darkblue] (0.28,-.3) to (-0.28,.4);
	\draw[->,thick,darkblue] (-0.28,-.3) to (0.28,.4);
\end{tikzpicture}
}
\:\:\:
\mathord{
\begin{tikzpicture}[baseline = -0.9mm]
	\draw[<-,thick,darkblue] (0.4,0.2) to[out=-90, in=0] (0.1,-.2);
	\draw[-,thick,darkblue] (0.1,-.2) to[out = 180, in = -90] (-0.2,0.2);
\end{tikzpicture}
}\:\right]
:\up\otimes \down \oplus \unit
\stackrel{\sim}{\rightarrow}
\down \otimes\up
$$
in $\mathcal H$, mirroring
the Mackey decomposition
$$
\ind_{n-1}^n \circ \res_{n-1}^n \oplus \id_n
\cong 
\res^{n+1}_n \circ \ind_{n}^{n+1}
$$
at the level of representation theory of the symmetric groups.
There have been several subsequent generalizations of Khovanov's work,
including a $q$-deformation \cite{LS}, a version of Heisenberg
category for wreath
product algebras associated to finite subgroups of $SL_2(\C)$
\cite{CL},
and an odd analog incorporating a Clifford superalgebra
\cite{HS}.

To explain the name ``Heisenberg category,''
let
$\mathfrak{h}$ be the infinite-dimensional Heisenberg algebra, i.e., the complex Lie algebra 
with basis $\{c, p_n, q_n\:|\:n \geq 1\}$ 
and multiplication given by
$$
[p_m, p_n] = [q_m, q_n] = [c,p_n]=[c,q_n] = 0,\qquad
[p_m, q_n] = \delta_{m,n} m c.
$$
Khovanov
constructed an algebra homomorphism 
from $U(\mathfrak{h})$ specialized at central charge $c=-1$ 
to
the complexified Grothendieck ring 
$\C\otimes_{\Z} K_0(\Kar(\mathcal H))$
of the additive
Karoubi envelope $\Kar(\mathcal H)$ of $\mathcal H$. He proved that his map is
injective, and conjectured that it is actually
an isomorphism. This conjecture is still open.

We remark also that the trace of Khovanov's category
and of its $q$-deformed version have recently been
computed; see \cite{CLLS, CLLSS}.

The group algebra of the symmetric group is the level one case of
a family of finite-dimensional algebras: the cyclotomic
quotients of degenerate affine Hecke algebras associated to
symmetric groups.
For cyclotomic quotients of level $\ell > 0$, the Mackey theorem instead
takes the form
$$
\ind_{n-1}^n \circ \res^{n}_{n-1} \oplus (\id_n)^{\oplus
  \ell} 
\cong \res^{n+1}_n \circ \ind_{n}^{n+1},
$$
e.g., see \cite[Theorem 7.6.2]{Kbook}.
Mackaay and Savage \cite{MS} have extended Khovanov's
construction to this setting, defining Heisenberg categories 
for all $\ell > 0$, with the case $\ell=1$ recovering Khovanov's
original category.
They also constructed an injective homomorphism from $U(\mathfrak{h})$
specialized at central charge $c=-\ell$ to 
the
complexified Grothendieck ring of the additive Karoubi envelope of their category,
and conjectured that this map is an
isomorphism. Again, this more general conjecture remains open.

In \cite{BCNR}, motivated by quite different considerations, the author jointly with
Comes, Nash and Reynolds introduced another diagrammatically-defined 
monoidal category we called the
{\em affine oriented Brauer category} $\AOB$; the
endomorphism algebras of objects in $\AOB$ are the {\em affine walled
  Brauer algebras} of \cite{RS}. In fact, the affine oriented Brauer
category is the Heisenberg category for central charge zero. 
To make this connection explicit, and also to streamline
the approach of Mackaay and Savage, we
propose here a simplified definition of Heisenberg
category for an arbitrary central charge $k \in \Z$.
Our new formulation is similar in spirit to Rouquier's definition of Kac-Moody
2-category from \cite{Rou} (as opposed to the Khovanov-Lauda definition
from \cite{KL3}); see also \cite{B}.

\begin{definition}\label{maindef}
Fix a commutative ground ring $\k$.
The {\em Heisenberg category} 
$\H_k$
of central charge $k \in \Z$ 
is the strict $\k$-linear monoidal category 
generated by objects
$\E$ and $\F$,
and
morphisms
$x:\E \rightarrow \E$,
$s:\E\otimes \E \rightarrow \E\otimes \E$, 
$c: \unit \rightarrow \F\otimes \E$ and
$d:\E \otimes \F \rightarrow \unit$ subject to certain relations.
To record these relations, we adopt the usual string calculus for strict
monoidal categories, representing the
generating morphisms by the diagrams
$$
x= 
\mathord{
\begin{tikzpicture}[baseline = 0]
	\draw[->,thick,darkblue] (0.08,-.3) to (0.08,.4);
      \node at (0.08,0.05) {$\dot$};
\end{tikzpicture}
}
\:,\qquad\quad
s = 
\mathord{
\begin{tikzpicture}[baseline = 0]
	\draw[->,thick,darkblue] (0.28,-.3) to (-0.28,.4);
	\draw[thick,darkblue,->] (-0.28,-.3) to (0.28,.4);
\end{tikzpicture}
}
\:,
\qquad\quad c = 
\mathord{
\begin{tikzpicture}[baseline = 1mm]
	\draw[<-,thick,darkblue] (0.4,0.4) to[out=-90, in=0] (0.1,0);
	\draw[-,thick,darkblue] (0.1,0) to[out = 180, in = -90] (-0.2,0.4);
\end{tikzpicture}
}
\:,
\qquad\quad
d = 
\mathord{
\begin{tikzpicture}[baseline = 1mm]
	\draw[<-,thick,darkblue] (0.4,0) to[out=90, in=0] (0.1,0.4);
	\draw[-,thick,darkblue] (0.1,0.4) to[out = 180, in = 90] (-0.2,0);
\end{tikzpicture}
}\:.
$$
The horizontal composition
$a \otimes b$ of two morphisms is
$a$ drawn to the left of $b$, and the vertical
composition $a \circ b$ is $a$ drawn above $b$ (assuming this makes sense).
We also denote the $n$th power $x^{\circ n}$ of $x$ under vertical composition 
diagrammatically by
labeling the dot with the multiplicity $n$, and define
$t:\E\otimes \F \rightarrow \F\otimes \E$ from
\begin{align}\label{ha}
t&=\mathord{
\begin{tikzpicture}[baseline = 0]
	\draw[<-,thick,darkblue] (0.28,-.3) to (-0.28,.4);
	\draw[->,thick,darkblue] (-0.28,-.3) to (0.28,.4);
\end{tikzpicture}
}
:=
\mathord{
\begin{tikzpicture}[baseline = 0]
	\draw[->,thick,darkblue] (0.3,-.5) to (-0.3,.5);
	\draw[-,thick,darkblue] (-0.2,-.2) to (0.2,.3);
        \draw[-,thick,darkblue] (0.2,.3) to[out=50,in=180] (0.5,.5);
        \draw[->,thick,darkblue] (0.5,.5) to[out=0,in=90] (0.8,-.5);
        \draw[-,thick,darkblue] (-0.2,-.2) to[out=230,in=0] (-0.5,-.5);
        \draw[-,thick,darkblue] (-0.5,-.5) to[out=180,in=-90] (-0.8,.5);
\end{tikzpicture}
}\,.
\end{align}
Then we impose three sets of relations: degenerate Hecke
relations, right adjunction relations, and
the inversion relation.
The degenerate Hecke relations are as follows\footnote{The final one
  of these relations is in parentheses to indicate that it is a
  consequence of the other relations; we have
  included it just for convenience.}:
\begin{align}\label{hecke}
\mathord{
\begin{tikzpicture}[baseline = -1mm]
	\draw[->,thick,darkblue] (0.28,0) to[out=90,in=-90] (-0.28,.6);
	\draw[->,thick,darkblue] (-0.28,0) to[out=90,in=-90] (0.28,.6);
	\draw[-,thick,darkblue] (0.28,-.6) to[out=90,in=-90] (-0.28,0);
	\draw[-,thick,darkblue] (-0.28,-.6) to[out=90,in=-90] (0.28,0);
\end{tikzpicture}
}=
\mathord{
\begin{tikzpicture}[baseline = -1mm]
	\draw[->,thick,darkblue] (0.18,-.6) to (0.18,.6);
	\draw[->,thick,darkblue] (-0.18,-.6) to (-0.18,.6);
\end{tikzpicture}
}\,,\qquad
\mathord{
\begin{tikzpicture}[baseline = -1mm]
	\draw[<-,thick,darkblue] (0.45,.6) to (-0.45,-.6);
	\draw[->,thick,darkblue] (0.45,-.6) to (-0.45,.6);
        \draw[-,thick,darkblue] (0,-.6) to[out=90,in=-90] (-.45,0);
        \draw[->,thick,darkblue] (-0.45,0) to[out=90,in=-90] (0,0.6);
\end{tikzpicture}
}
&=
\mathord{
\begin{tikzpicture}[baseline = -1mm]
	\draw[<-,thick,darkblue] (0.45,.6) to (-0.45,-.6);
	\draw[->,thick,darkblue] (0.45,-.6) to (-0.45,.6);
        \draw[-,thick,darkblue] (0,-.6) to[out=90,in=-90] (.45,0);
        \draw[->,thick,darkblue] (0.45,0) to[out=90,in=-90] (0,0.6);
\end{tikzpicture}
}\,,\quad
&\mathord{
\begin{tikzpicture}[baseline = -1mm]
	\draw[<-,thick,darkblue] (0.25,.4) to (-0.25,-.4);
	\draw[->,thick,darkblue] (0.25,-.4) to (-0.25,.4);
     \node at (-0.12,0.2) {$\dot$};
\end{tikzpicture}
}
-
\mathord{
\begin{tikzpicture}[baseline = -1mm]
	\draw[<-,thick,darkblue] (0.25,.4) to (-0.25,-.4);
	\draw[->,thick,darkblue] (0.25,-.4) to (-0.25,.4);
     \node at (0.12,-0.2) {$\dot$};
\end{tikzpicture}
}
&=\mathord{
\begin{tikzpicture}[baseline = -1mm]
 	\draw[->,thick,darkblue] (0.08,-.4) to (0.08,.4);
	\draw[->,thick,darkblue] (-0.28,-.4) to (-0.28,.4);
\end{tikzpicture}
}\:\left(=\mathord{
\begin{tikzpicture}[baseline = -1mm]
	\draw[<-,thick,darkblue] (0.25,.4) to (-0.25,-.4);
	\draw[->,thick,darkblue] (0.25,-.4) to (-0.25,.4);
     \node at (-0.12,-0.2) {$\dot$};
\end{tikzpicture}
}
-
\mathord{
\begin{tikzpicture}[baseline = -1mm]
	\draw[<-,thick,darkblue] (0.25,.4) to (-0.25,-.4);
	\draw[->,thick,darkblue] (0.25,-.4) to (-0.25,.4);
     \node at (0.12,0.2) {$\dot$};
\end{tikzpicture}
}\right)\,.
\end{align}
The right adjunction relations say that
\begin{align}\label{rightadj}
\mathord{
\begin{tikzpicture}[baseline = 0]
  \draw[->,thick,darkblue] (0.3,0) to (0.3,.4);
	\draw[-,thick,darkblue] (0.3,0) to[out=-90, in=0] (0.1,-0.4);
	\draw[-,thick,darkblue] (0.1,-0.4) to[out = 180, in = -90] (-0.1,0);
	\draw[-,thick,darkblue] (-0.1,0) to[out=90, in=0] (-0.3,0.4);
	\draw[-,thick,darkblue] (-0.3,0.4) to[out = 180, in =90] (-0.5,0);
  \draw[-,thick,darkblue] (-0.5,0) to (-0.5,-.4);
\end{tikzpicture}
}
&=
\mathord{\begin{tikzpicture}[baseline=0]
  \draw[->,thick,darkblue] (0,-0.4) to (0,.4);
\end{tikzpicture}
}\,,
&\mathord{
\begin{tikzpicture}[baseline = 0]
  \draw[->,thick,darkblue] (0.3,0) to (0.3,-.4);
	\draw[-,thick,darkblue] (0.3,0) to[out=90, in=0] (0.1,0.4);
	\draw[-,thick,darkblue] (0.1,0.4) to[out = 180, in = 90] (-0.1,0);
	\draw[-,thick,darkblue] (-0.1,0) to[out=-90, in=0] (-0.3,-0.4);
	\draw[-,thick,darkblue] (-0.3,-0.4) to[out = 180, in =-90] (-0.5,0);
  \draw[-,thick,darkblue] (-0.5,0) to (-0.5,.4);
\end{tikzpicture}
}
&=
\mathord{\begin{tikzpicture}[baseline=0]
  \draw[<-,thick,darkblue] (0,-0.4) to (0,.4);
\end{tikzpicture}
}\,.
\end{align}
Finally, the inversion relation asserts
that the following matrix of
morphisms
is an isomorphism
in the additive envelope of $\H_k$:
\begin{align}
\label{invrel1}
\left[\:
\mathord{
\begin{tikzpicture}[baseline = 0]
	\draw[<-,thick,darkblue] (0.28,-.3) to (-0.28,.4);
	\draw[->,thick,darkblue] (-0.28,-.3) to (0.28,.4);
   \end{tikzpicture}
}
\:\:\:
\mathord{
\begin{tikzpicture}[baseline = 1mm]
	\draw[<-,thick,darkblue] (0.4,0) to[out=90, in=0] (0.1,0.4);
	\draw[-,thick,darkblue] (0.1,0.4) to[out = 180, in = 90] (-0.2,0);
\end{tikzpicture}
}\:\:\:\mathord{
\begin{tikzpicture}[baseline = 1mm]
	\draw[<-,thick,darkblue] (0.4,0) to[out=90, in=0] (0.1,0.4);
	\draw[-,thick,darkblue] (0.1,0.4) to[out = 180, in = 90] (-0.2,0);
      \node at (-0.15,0.2) {$\dot$};
\end{tikzpicture}
}\:\:\:\cdots
\:\:\:
\mathord{
\begin{tikzpicture}[baseline = 1mm]
	\draw[<-,thick,darkblue] (0.4,0) to[out=90, in=0] (0.1,0.4);
	\draw[-,thick,darkblue] (0.1,0.4) to[out = 180, in = 90] (-0.2,0);
     \node at (-0.55,0.2) {$\color{darkblue}\scriptstyle{k-1}$};
      \node at (-0.15,0.2) {$\dot$};
\end{tikzpicture}
}\:\right]^T
&:
\E \otimes \F\stackrel{\sim}{\rightarrow}
\F \otimes \E \oplus \unit^{\oplus k}
&&\text{if $k \geq
  0$},\\
\left[\:
\mathord{
\begin{tikzpicture}[baseline = 0]
	\draw[<-,thick,darkblue] (0.28,-.3) to (-0.28,.4);
	\draw[->,thick,darkblue] (-0.28,-.3) to (0.28,.4);
\end{tikzpicture}
}\:\:\:
\mathord{
\begin{tikzpicture}[baseline = -0.9mm]
	\draw[<-,thick,darkblue] (0.4,0.2) to[out=-90, in=0] (0.1,-.2);
	\draw[-,thick,darkblue] (0.1,-.2) to[out = 180, in = -90] (-0.2,0.2);
\end{tikzpicture}
}
\:\:\:
\mathord{
\begin{tikzpicture}[baseline = -0.9mm]
	\draw[<-,thick,darkblue] (0.4,0.2) to[out=-90, in=0] (0.1,-.2);
	\draw[-,thick,darkblue] (0.1,-.2) to[out = 180, in = -90] (-0.2,0.2);
      \node at (0.38,0) {$\dot$};
\end{tikzpicture}
}
\:\:\:\cdots
\:\:\:
\mathord{
\begin{tikzpicture}[baseline = -0.9mm]
	\draw[<-,thick,darkblue] (0.4,0.2) to[out=-90, in=0] (0.1,-.2);
	\draw[-,thick,darkblue] (0.1,-.2) to[out = 180, in = -90] (-0.2,0.2);
     \node at (0.83,0) {$\color{darkblue}\scriptstyle{-k-1}$};
      \node at (0.38,0) {$\dot$};
\end{tikzpicture}
}
\right]
&:\E \otimes \F \oplus 
\unit^{\oplus (-k)}
\stackrel{\sim}{\rightarrow}
 \F \otimes  \E&&\text{if $k <
  0$}.
\label{invrel2}
\end{align}
\end{definition}


In the special case $k=0$, the inversion
relation
means that one should adjoin another generating morphism
$t':\F\otimes \E \rightarrow \E\otimes \F$, 
represented by
$$t'=
\mathord{
\begin{tikzpicture}[baseline = 0]
	\draw[->,thick,darkblue] (0.28,-.3) to (-0.28,.4);
	\draw[<-,thick,darkblue] (-0.28,-.3) to (0.28,.4);
\end{tikzpicture}
},
$$
subject to the following relations asserting that $t'$ is a two-sided
inverse to $t$:
\begin{align*}
\mathord{
\begin{tikzpicture}[baseline = 0mm]
	\draw[->,thick,darkblue] (0.28,0) to[out=90,in=-90] (-0.28,.6);
	\draw[-,thick,darkblue] (-0.28,0) to[out=90,in=-90] (0.28,.6);
	\draw[<-,thick,darkblue] (0.28,-.6) to[out=90,in=-90] (-0.28,0);
	\draw[-,thick,darkblue] (-0.28,-.6) to[out=90,in=-90] (0.28,0);
\end{tikzpicture}
}
&=\mathord{
\begin{tikzpicture}[baseline = 0]
	\draw[<-,thick,darkblue] (0.08,-.6) to (0.08,.6);
	\draw[->,thick,darkblue] (-0.28,-.6) to (-0.28,.6);
\end{tikzpicture}
}\:,
&\mathord{
\begin{tikzpicture}[baseline = 0mm]
	\draw[-,thick,darkblue] (0.28,0) to[out=90,in=-90] (-0.28,.6);
	\draw[->,thick,darkblue] (-0.28,0) to[out=90,in=-90] (0.28,.6);
	\draw[-,thick,darkblue] (0.28,-.6) to[out=90,in=-90] (-0.28,0);
	\draw[<-,thick,darkblue] (-0.28,-.6) to[out=90,in=-90] (0.28,0);
\end{tikzpicture}
}
&=
\mathord{
\begin{tikzpicture}[baseline = 0]
	\draw[->,thick,darkblue] (0.08,-.6) to (0.08,.6);
	\draw[<-,thick,darkblue] (-0.28,-.6) to (-0.28,.6);
\end{tikzpicture}
}\:.
\end{align*}
Up to reflecting diagrams in a vertical axis, this is exactly the
definition of the {affine oriented Brauer category}
$\mathcal{AOB}$
from \cite{BCNR}. Thus, there is a monoidal isomorphism
$\H_0
\cong\mathcal{AOB}^{\operatorname{rev}}$.

When $k \neq 0$, the inversion relation appearing in Definition~\ref{maindef}
is much harder to interpret. We will analyze it systematically in
the main part of this article. We summarize the situation with
the following two theorems.

\begin{theorem}\label{thm1}
There are unique 
morphisms $c':\unit \rightarrow \E\otimes \F$ and $d':\F\otimes \E \rightarrow
\unit$ in $\H_k$,
drawn as
$$
c' = 
\mathord{

}\,,
$$
such that the following relations hold:
\begin{align}\label{pos}
\mathord{
%
}= 
\delta_{r,-k-1} 1_\unit
\text{ if $0 \leq r < -k$.}
\end{align}
Moreover,
$\H_k$ can be presented equivalently as the strict $\k$-linear monoidal category
generated by the objects $\E,\F$ and
morphisms
$x,s,c,d,c',d'$ subject only to the relations
(\ref{hecke})--(\ref{rightadj}) and (\ref{pos})--(\ref{leftcurl}).
In these relations, as well as
the rightward crossing $t$ defined by (\ref{ha}), we
have used the
leftward crossing $t':\F\otimes \E \rightarrow \E\otimes \F$
 defined by
\begin{equation}\label{tprime}
t'=
\mathord{
%
}\,,\end{equation}
and the negatively dotted bubbles
defined by
\begin{align}\label{d1}
\mathord{
%
}\right)_{i,j=1,\dots,r}
&&\text{if $r \leq -k$,}\label{d2}
\end{align}
interpreting the determinants as $1$ if $r=0$ and $0$ if $r
< 0$.
\end{theorem}

\begin{theorem}\label{thm2}
Using the notation from Theorem~\ref{thm1},
the following relations are consequences of the defining relations.
\begin{itemize}
\item
[(i)] (``Infinite Grassmannian relations'') 
For all $r \in \Z\,$:
\begin{align}\label{fed1}
\mathord{
%
}\,.
\end{align}
\item[(vi)] (``Alternating braid relation'')
\begin{align}
\mathord{

}
&=
\left\{
\begin{array}{cl}
\displaystyle\sum_{r,s,t \geq 0}
\mathord{
\begin{tikzpicture}[baseline = 0]
	\draw[-,thick,darkblue] (0.3,0.6) to[out=-90, in=0] (0,0.2);
	\draw[->,thick,darkblue] (0,0.2) to[out = 180, in = -90] (-0.3,0.6);
  \draw[->,thick,darkblue] (-0.8,0) to[out=90,in=0] (-1,0.2);
  \draw[-,thick,darkblue] (-1,0.2) to[out=180,in=90] (-1.2,0);
\draw[-,thick,darkblue] (-1.2,0) to[out=-90,in=180] (-1,-0.2);
  \draw[-,thick,darkblue] (-1,-0.2) to[out=0,in=-90] (-0.8,0);
   \node at (-0.8,0) {$\dot$};
   \node at (-.14,0) {$\color{darkblue}\scriptstyle{-r-s-t-3}$};
   \node at (-0.23,0.33) {$\dot$};
   \node at (-0.43,0.33) {$\color{darkblue}\scriptstyle{r}$};
	\draw[<-,thick,darkblue] (0.3,-.6) to[out=90, in=0] (0,-0.2);
	\draw[-,thick,darkblue] (0,-0.2) to[out = 180, in = 90] (-0.3,-.6);
   \node at (-0.25,-0.4) {$\dot$};
   \node at (-.4,-.4) {$\color{darkblue}\scriptstyle{s}$};
	\draw[->,thick,darkblue] (.98,-0.6) to (.98,0.6);
   \node at (0.98,0) {$\dot$};
   \node at (1.15,0) {$\color{darkblue}\scriptstyle{t}$};
\end{tikzpicture}
}
&\text{if $k \geq 2$,}\\\\
0&\text{if $-1 \leq k \leq 1$,}\\\\
\end{array}\right.
\label{altbraid}
\end{align}
\end{itemize}
\end{theorem}

Part (ii) of Theorem~\ref{thm2} implies
that the monoidal category $\H_k$ is {\em rigid}, i.e., any object $X$ has
both 
a right dual $X^*$ (with its structure maps
$X \otimes X^* \rightarrow \unit \rightarrow X^* \otimes X$) 
and
a left dual ${^*\!}X$ (with its structure maps
${^*\!}X\otimes X  \rightarrow \unit \rightarrow X \otimes {^*\!}X$).
In fact, there is a canonical
choice for both duals, by attaching the appropriately oriented cups and caps 
as indicated below:
$$
\mathord{
%
}.
$$
Then part (iii) of the theorem
shows that the right and left mates of $x$ are equal, as are the right
and left mates of $s$.
We denote these by
$x':\F \rightarrow \F$ and $s':\F\otimes \F \rightarrow \F\otimes \F$,
respectively,
and represent them diagrammatically by
\begin{align*}
x' &= \mathord{
%
}\:.
\end{align*}
It follows that the functors $(-)^*$ and ${^*}(-)$ defined by taking right and
left duals/mates in the canonical way actually coincide;
they are both defined by rotating diagrams through $180^\circ$.
Thus, we have equipped $\H_k$ with a {\em strictly pivotal structure}.

Now we can explain the relationship between the category $\H_k$
and the Heisenberg categories 
already appearing in the literature.
By a special case of the bubble slide relations in Theorem~\ref{thm2}(v),
the
lowest degree bubble 
$\mathord{%
}
$ is {\em strictly central} in the sense that
\begin{align*}
\mathord{
%
}
\,.
\end{align*}
This means that it is natural to specialize 
$\mathord{%
}$ to some scalar $\delta \in \k$.
We denote the resulting monoidal category by $\H_k(\delta)$.

\begin{theorem}\label{thm3}
The 
Heisenberg category $\tilde{\mathcal H}^\lambda$
defined by Mackaay and Savage in \cite{MS} 
is isomorphic to the additive envelope of 
$\H_{k}(\delta)$, taking
$k:=-\sum_i \lambda_i$
and $\delta := \sum_{i} i \lambda_i$.
In particular,
the
Heisenberg category $\mathcal H$ introduced originally by Khovanov
in \cite{K} is isomorphic to
the additive envelope of $\H_{-1}(0)$.
\end{theorem}

\begin{remark}
The above results give {\em two} new presentations for Khovanov's
Heisenberg category $\mathcal H$, i.e., the additive envelope of our $\H_{-1}(0)$:
\begin{enumerate}
\item
The first presentation, which is essentially
Definition~\ref{maindef}, asserts that $\mathcal H$ is the strict
additive $\k$-linear monoidal
category generated by objects
$\up$ and $\down$ and the
morphisms $x,s,c$ and $d$, subject to the relations (\ref{hecke}) and
(\ref{rightadj}), the relation
(\ref{invrel2}) asserting
that
$$
\left[\:
\mathord{
\begin{tikzpicture}[baseline = 0]
	\draw[<-,thick,darkblue] (0.28,-.3) to (-0.28,.4);
	\draw[->,thick,darkblue] (-0.28,-.3) to (0.28,.4);
\end{tikzpicture}
}\:\:\:
\mathord{
\begin{tikzpicture}[baseline = -0.9mm]
	\draw[<-,thick,darkblue] (0.4,0.2) to[out=-90, in=0] (0.1,-.2);
	\draw[-,thick,darkblue] (0.1,-.2) to[out = 180, in = -90] (-0.2,0.2);
\end{tikzpicture}
}
\:
\right]
:\E \otimes \F \oplus 
\unit
\stackrel{\sim}{\rightarrow}
 \F \otimes  \E
$$
is an isomorphism where the rightward crossing is defined by (\ref{ha}), and the relation
$$\mathord{\begin{tikzpicture}[baseline = -1mm]
  \draw[-,thick,darkblue] (0,0.2) to[out=180,in=90] (-.2,0);
  \draw[->,thick,darkblue] (0.2,0) to[out=90,in=0] (0,.2);
 \draw[-,thick,darkblue] (-.2,0) to[out=-90,in=180] (0,-0.2);
  \draw[-,thick,darkblue] (0,-0.2) to[out=0,in=-90] (0.2,0);
      \node at (0.2,0) {$\dot$};
\end{tikzpicture}
} = 0$$
where the leftward cap is defined
from
$$
\left[\begin{array}{l}
\displaystyle\mathord{
\begin{tikzpicture}[baseline = 0]
	\draw[->,thick,darkblue] (0.28,-.3) to (-0.28,.4);
	\draw[<-,thick,darkblue] (-0.28,-.3) to (0.28,.4);
\end{tikzpicture}
}\\\\
\mathord{
\begin{tikzpicture}[baseline = 1mm]
	\draw[-,thick,darkblue] (0.4,0) to[out=90, in=0] (0.1,0.4);
	\draw[->,thick,darkblue] (0.1,0.4) to[out = 180, in = 90] (-0.2,0);
\end{tikzpicture}
}
\end{array}\right] := \left[\:
\mathord{
\begin{tikzpicture}[baseline = 0]
	\draw[<-,thick,darkblue] (0.28,-.3) to (-0.28,.4);
	\draw[->,thick,darkblue] (-0.28,-.3) to (0.28,.4);
\end{tikzpicture}
}\:\:\:
\mathord{
\begin{tikzpicture}[baseline = -0.9mm]
	\draw[<-,thick,darkblue] (0.4,0.2) to[out=-90, in=0] (0.1,-.2);
	\draw[-,thick,darkblue] (0.1,-.2) to[out = 180, in = -90] (-0.2,0.2);
\end{tikzpicture}
}
\:
\right]
^{-1}.
$$
\vspace{-2mm}

\noindent
The leftward cup may also be recovered from
$\mathord{
\begin{tikzpicture}[baseline = .5mm]
	\draw[-,thick,darkblue] (0.3,0.3) to[out=-90, in=0] (0.1,0);
	\draw[->,thick,darkblue] (0.1,0) to[out = 180, in = -90] (-0.1,0.3);
\end{tikzpicture}
}:=
\mathord{
\begin{tikzpicture}[baseline = -1mm]
	\draw[-,thick,darkblue] (0.2,.4) to[out=240,in=90] (-0.2,-.1);
	\draw[<-,thick,darkblue] (-0.2,.4) to[out=300,in=90] (0.2,-0.1);
	\draw[-,thick,darkblue] (0.2,-0.1) to[out=-90, in=0] (0,-0.3);
	\draw[-,thick,darkblue] (0,-0.3) to[out = 180, in = -90] (-0.2,-0.1);
      \node at (0.2,-0.1) {$\dot$};
\end{tikzpicture}
}\:.$
\item
The second presentation, which is a simplification of the
presentation from Theorem~\ref{thm2},
asserts that $\mathcal H$ is generated by objects $\up$ and
$\down$
and the morphisms $s,c,d,c',d'$ subject to the first two relations from
(\ref{hecke}),
the relations (\ref{rightadj}), and four additional relations:
\begin{align*}
\mathord{
\begin{tikzpicture}[baseline = 0mm]
	\draw[->,thick,darkblue] (0.28,0) to[out=90,in=-90] (-0.28,.6);
	\draw[-,thick,darkblue] (-0.28,0) to[out=90,in=-90] (0.28,.6);
	\draw[<-,thick,darkblue] (0.28,-.6) to[out=90,in=-90] (-0.28,0);
	\draw[-,thick,darkblue] (-0.28,-.6) to[out=90,in=-90] (0.28,0);
\end{tikzpicture}
}
&=\mathord{
\begin{tikzpicture}[baseline = 0]
	\draw[<-,thick,darkblue] (0.08,-.6) to (0.08,.6);
	\draw[->,thick,darkblue] (-0.28,-.6) to (-0.28,.6);
\end{tikzpicture}
}\:,
&
\mathord{
\begin{tikzpicture}[baseline = 0mm]
	\draw[-,thick,darkblue] (0.28,0) to[out=90,in=-90] (-0.28,.6);
	\draw[->,thick,darkblue] (-0.28,0) to[out=90,in=-90] (0.28,.6);
	\draw[-,thick,darkblue] (0.28,-.6) to[out=90,in=-90] (-0.28,0);
	\draw[<-,thick,darkblue] (-0.28,-.6) to[out=90,in=-90] (0.28,0);
\end{tikzpicture}
}
&=\mathord{
\begin{tikzpicture}[baseline = 0]
	\draw[->,thick,darkblue] (0.08,-.6) to (0.08,.6);
	\draw[<-,thick,darkblue] (-0.28,-.6) to (-0.28,.6);
\end{tikzpicture}
}
-\mathord{
\begin{tikzpicture}[baseline=-.5mm]
	\draw[<-,thick,darkblue] (0.3,0.6) to[out=-90, in=0] (0,.1);
	\draw[-,thick,darkblue] (0,.1) to[out = 180, in = -90] (-0.3,0.6);
	\draw[-,thick,darkblue] (0.3,-.6) to[out=90, in=0] (0,-0.1);
	\draw[->,thick,darkblue] (0,-0.1) to[out = 180, in = 90] (-0.3,-.6);
\end{tikzpicture}}\,,
&
\mathord{
\begin{tikzpicture}[baseline = -0.5mm]
	\draw[<-,thick,darkblue] (0,0.6) to (0,0.3);
	\draw[-,thick,darkblue] (0,0.3) to [out=-90,in=0] (-.3,-0.2);
	\draw[-,thick,darkblue] (-0.3,-0.2) to [out=180,in=-90](-.5,0);
	\draw[-,thick,darkblue] (-0.5,0) to [out=90,in=180](-.3,0.2);
	\draw[-,thick,darkblue] (-0.3,.2) to [out=0,in=90](0,-0.3);
	\draw[-,thick,darkblue] (0,-0.3) to (0,-0.6);
\end{tikzpicture}
}&=
0,
&
\mathord{
\begin{tikzpicture}[baseline = 1.25mm]
  \draw[->,thick,darkblue] (0.2,0.2) to[out=90,in=0] (0,.4);
  \draw[-,thick,darkblue] (0,0.4) to[out=180,in=90] (-.2,0.2);
\draw[-,thick,darkblue] (-.2,0.2) to[out=-90,in=180] (0,0);
  \draw[-,thick,darkblue] (0,0) to[out=0,in=-90] (0.2,0.2);
\end{tikzpicture}
}&= 1_\unit.
\end{align*}
The rightward and leftward crossings used here are shorthands for the
morphisms defined by (\ref{ha}) and (\ref{tprime}), respectively.
Then $x$ may be defined from
$\mathord{
\begin{tikzpicture}[baseline = -.5mm]
	\draw[->,thick,darkblue] (0.08,-.3) to (0.08,.3);
      \node at (0.08,0) {$\dot$};
\end{tikzpicture}
}:=
\mathord{
\begin{tikzpicture}[baseline = -0.5mm]
	\draw[<-,thick,darkblue] (0,0.3) to (0,0.2);
	\draw[-,thick,darkblue] (0,0.2) to [out=-90,in=180] (.2,-0.15);
	\draw[-,thick,darkblue] (0.2,-0.15) to [out=0,in=-90](.35,0);
	\draw[-,thick,darkblue] (0.35,0) to [out=90,in=0](.2,0.15);
	\draw[-,thick,darkblue] (0.2,.15) to [out=180,in=90](0,-0.2);
	\draw[-,thick,darkblue] (0,-0.2) to (0,-0.3);
\end{tikzpicture}
}\,$; the third relation
from (\ref{hecke}) holds automatically.
\end{enumerate}
The presentation (2) is almost the same as
Khovanov's original definition. Khovanov's formulation also implicitly incorporated some additional cyclicity relations, which our results show are redundant, i.e., they are implied by the other relations.
\end{remark}

Let $\Sym$ be the algebra of symmetric functions.
Recall this is an infinite rank polynomial algebra 
generated freely by either the complete symmetric functions
$\{h_r\}_{r \geq 1}$ or the elementary symmetric functions
$\{e_r\}_{r \geq 1}$; we also let $h_0 = e_0 = 1$ and interpret
$h_r$ and $e_r$ as $0$ when $r < 0$.
Let
\begin{equation}\label{savvy}
\beta:\Sym \rightarrow \End_{\H_k}(\unit)
\end{equation}
be the algebra homomorphism defined by
declaring that 
$$
\left\{\begin{array}{rl}
\beta(e_r) := 
-\mathord{
\begin{tikzpicture}[baseline = 1.25mm]
  \draw[<-,thick,darkblue] (0,0.4) to[out=180,in=90] (-.2,0.2);
  \draw[-,thick,darkblue] (0.2,0.2) to[out=90,in=0] (0,.4);
 \draw[-,thick,darkblue] (-.2,0.2) to[out=-90,in=180] (0,0);
  \draw[-,thick,darkblue] (0,0) to[out=0,in=-90] (0.2,0.2);
   \node at (-0.2,0.2) {$\dot$};
   \node at (-0.7,0.2) {$\color{darkblue}\scriptstyle{r+k-1}$};
\end{tikzpicture}}\phantom{(-1)^r}&\text{if $k \geq 0$,}\\
\beta(h_r) :=(-1)^r \mathord{
\begin{tikzpicture}[baseline = 1.25mm]
  \draw[->,thick,darkblue] (0.2,0.2) to[out=90,in=0] (0,.4);
  \draw[-,thick,darkblue] (0,0.4) to[out=180,in=90] (-.2,0.2);
\draw[-,thick,darkblue] (-.2,0.2) to[out=-90,in=180] (0,0);
  \draw[-,thick,darkblue] (0,0) to[out=0,in=-90] (0.2,0.2);
   \node at (0.2,0.2) {$\dot$};
   \node at (0.7,0.2) {$\color{darkblue}\scriptstyle{r-k-1}$};
\end{tikzpicture}
}\phantom{-}&\text{if $k < 0$.}
\end{array}\right.
$$
Then the relations from Theorem~\ref{thm2}(i)
imply that 
$$
\left\{\begin{array}{rl}
\beta(h_r) = 
(-1)^r \mathord{
\begin{tikzpicture}[baseline = 1.25mm]
  \draw[->,thick,darkblue] (0.2,0.2) to[out=90,in=0] (0,.4);
  \draw[-,thick,darkblue] (0,0.4) to[out=180,in=90] (-.2,0.2);
\draw[-,thick,darkblue] (-.2,0.2) to[out=-90,in=180] (0,0);
  \draw[-,thick,darkblue] (0,0) to[out=0,in=-90] (0.2,0.2);
   \node at (0.2,0.2) {$\dot$};
   \node at (0.7,0.2) {$\color{darkblue}\scriptstyle{r-k-1}$};
\end{tikzpicture}
}\phantom{-}&\text{if $k \geq 0$,}\\
\beta(e_r) = -\mathord{
\begin{tikzpicture}[baseline = 1.25mm]
  \draw[<-,thick,darkblue] (0,0.4) to[out=180,in=90] (-.2,0.2);
  \draw[-,thick,darkblue] (0.2,0.2) to[out=90,in=0] (0,.4);
 \draw[-,thick,darkblue] (-.2,0.2) to[out=-90,in=180] (0,0);
  \draw[-,thick,darkblue] (0,0) to[out=0,in=-90] (0.2,0.2);
   \node at (-0.2,0.2) {$\dot$};
   \node at (-0.7,0.2) {$\color{darkblue}\scriptstyle{r+k-1}$};
\end{tikzpicture}
}\phantom{(-1)^r}&\text{if $k < 0$.}
\end{array}\right.
$$
In fact, $\beta$ is an isomorphism. This assertion is a consequence of the
{\em basis theorem}
for morphism spaces in $\H_k$, which we explain next.

Let $X = X_1\otimes \cdots \otimes X_r$ and $Y = Y_1 \otimes\cdots\otimes Y_s$
be two words in the letters $\up$ and $\down$, representing two
objects of $\H_k$.
By an {\em $(X,Y)$-matching}, we mean a 
bijection
$$
\{i\:|\:X_i = \up\}\sqcup\{j\:|\:Y_j = \down\}
\stackrel{\sim}{\rightarrow}
\{i\:|\:X_i = \down\}\sqcup\{j\:|\:Y_j = \up\}.
$$
By a {\em reduced lift} of an $(X,Y)$-matching, we mean a diagram
representing a morphism 
$X \rightarrow Y$ in $\H_k$ such that 
\begin{itemize}
\item the endpoints of each strand in the diagram are paired under the matching;
\item any two strands intersect at most once;
\item there are no self-intersections;
\item there are no dots or bubbles;
\item each strand has at most one critical point coming from a cup or
  cap.
\end{itemize}
Let $B(X,Y)$ be a set consisting of a reduced lift for each of the
$(X,Y)$-matchings.
For each element of
$B(X,Y)$, pick a distinguished point on each of its strands that is
away from crossings and critical points. Then let
$B_{\infty,\infty}(X,Y)$
be the set of all morphisms $\theta:X \rightarrow Y$ obtained from the
elements of $B(X,Y)$ by adding zero or more dots to each
strand at these distinguished points.

\begin{theorem}\label{bt}
For any $k \in \Z$ and objects
$X, Y$ in $\H_k$,
the space $\Hom_{\H_k}(X,Y)$ is a free right
$\Sym$-module with basis
$B_{\infty,\infty}(X,Y)$.
Here, the right action of $\Sym$ on morphisms is by $\theta \cdot p := \theta
\otimes \beta(p)$
for $\theta:X \rightarrow Y$ and $p \in \Sym$.
\end{theorem}

In particular, this implies that $\H_k \cong
\H_k(\delta)\otimes_\k \k[z]$ where $z$ denotes the bubble 
$\mathord{\begin{tikzpicture}[baseline = -1mm]
  \draw[-,thick,darkblue] (0,0.2) to[out=180,in=90] (-.2,0);
  \draw[-,thick,darkblue] (0.2,0) to[out=90,in=0] (0,.2);
 \draw[-,thick,darkblue] (-.2,0) to[out=-90,in=180] (0,-0.2);
  \draw[-,thick,darkblue] (0,-0.2) to[out=0,in=-90] (0.2,0);
\end{tikzpicture}
}$. 
It follows that $K_0(\H_k) \cong K_0(\H_k(\delta))$ for any $\delta
\in \k$. Combining this observation with Theorem~\ref{thm3}, we 
then restate
\cite[Theorem 4.4]{MS} as follows:
whenever $k \neq 0$
there is an algebra embedding 
\begin{equation}
U(\mathfrak{h}) / \langle c-k\rangle \hookrightarrow \C\otimes_{\Z}
K_0(\Kar(\H_k)).
\end{equation} 
As we mentioned already above, this embedding is conjectured to be
an isomorphism. There should be similar results when
$k=0$ too.

Theorem~\ref{bt} was proved already in case $k=0$ in \cite[Theorem 1.2]{BCNR}, 
by an argument based on the existence of 
a certain monoidal functor
from $\H_0$ to the category of $\k$-linear endofunctors of the
category of modules over the Lie algebra $\mathfrak{gl}_n(\k)$.
When $k \neq 0$, 
the theorem will instead be deduced from the basis
theorems proved in \cite[Proposition 5]{K} and 
\cite[Proposition 2.16]{MS}. The proofs 
in \cite{K, MS} depend crucially on the action of $\H_k$ on the category
of modules over the degenerate cyclotomic Hecke algebras mentioned earlier.
Since it highlights the usefulness of Definition~\ref{maindef}, we
give a self-contained construction of this action in the next paragraph.

Fix a monic polynomial $f(u)= u^\ell + z_1
u^{\ell-1}+\cdots+z_\ell\in \k[u]$ of degree $\ell > 0$ and set
$k := -\ell$.
Let $H_n$ be the degenerate affine Hecke algebra, that is, 
the tensor product $\k S_n \otimes \k[x_1,\dots,x_n]$
of the group algebra of the symmetric group with a polynomial algebra.
Multiplication in $H_n$ is defined so that $\k S_n$ and $\k[x_1,\dots,x_n]$
are subalgebras, and also
$$
x_{i+1} s_i = s_i x_i + 1, 
\qquad
x_i s_j = s_j x_i\:(i \neq j,j+1),
$$
where $s_j$ denotes the basic transposition $(j\:j\!+\!1)$.
Let $H_n^f$ be the
quotient of $H_n$ by the two-sided ideal generated by $f(x_1)$.
There is a natural embedding $H_n^f \hookrightarrow H_{n+1}^f$
sending $x_i,s_j \in H_n^f$ to the same elements of $H_{n+1}^f$.
Let 
\begin{align*}
\ind_n^{n+1} := H_{n+1}^f \otimes_{H_n^f} ?:H_n^f \Mod &\rightarrow
H_{n+1}^f \Mod,\\
\res_n^{n+1}: H_{n+1}^f\Mod &\rightarrow H_n^f \Mod
\end{align*} 
be the corresponding induction
and restriction functors.
The key assertion established in \cite{K, MS} is that there is a strict $\k$-linear monoidal functor
\begin{equation}\label{wolf}
\Psi_{f}:
\H_{k} \rightarrow \mathcal{E}\!nd_\k\left(\bigoplus_{n \geq 0} H_{n}^f\Mod\right)
\end{equation}
sending 
$\E$ (respectively, $\F$) to the $\k$-linear endofunctor that takes an $H_n^f$-module $M$ to
the $H_{n+1}^f$-module $\ind_n^{n+1} M$ (respectively,
to the $H_{n-1}^f$-module $\res^{n}_{n-1} M$, interpreted as zero in
case $n=0$).
On generating morphisms, $\Psi_f(x), \Psi_f(s), \Psi_f(c)$ and $\Psi_f(d)$ are
the natural transformations
defined on an $H_n^f$-module $M$ as follows:
\begin{itemize}
\item
$\Psi_f(x)_M:\ind_n^{n+1} M \rightarrow \ind_n^{n+1} M, \:
h \otimes m \mapsto h x_{n+1} \otimes m$;
\item
$\Psi_f(s)_M:\ind_n^{n+2} M\rightarrow \ind_n^{n+2} M, \:h\otimes m \mapsto h s_{n+1} \otimes m$, where we have identified
$\ind_{n+1}^{n+2} \circ \ind_n^{n+1}$ with $\ind_n^{n+2} := H_{n+2}^f
\otimes_{H_n^f} ?$ in the obvious way;
\item $\Psi_f(c)_M:M \rightarrow \res^{n+1}_n \circ \ind_n^{n+1} M, \:m
  \mapsto 1 \otimes m$;
\item $\Psi_f(d)_M:\ind_{n-1}^n \circ \res_{n-1}^n M \rightarrow M,
  \:h \otimes m
\mapsto hm$.
\end{itemize}
To prove this in our setting, we need to verify the three
sets of relations from Definition~\ref{maindef}. The first two are
almost immediate. For the
inversion relation, one calculates $\Psi_f(t)_M$
explicitly to see that it comes from the $(H_n^f, H_n^f)$-bimodule
homomorphism
$H_n^f \otimes_{H_{n-1}^f} H_n^f \rightarrow H_{n+1}^f, a \otimes b
\mapsto a s_n b$.
Thus, it suffices to show that 
the $(H_n^f,
H_n^f)$-bimodule homomorphism
\begin{align}\label{once}
H_n^f \otimes_{H_{n-1}^f} H_n^f \oplus \bigoplus_{r=0}^{\ell-1} H_n^f
&\rightarrow H_{n+1}^f,\\
(a\otimes b, c_0,c_1,\dots,c_{\ell-1})
&\mapsto a s_n b + \sum_{r=0}^{\ell-1} c_r x_{n+1}^r \notag
\end{align}
is an isomorphism, which is checked in the proof of
\cite[Lemma 7.6.1]{Kbook}. 
We remark further that $\Psi_f$ maps the bubble
$\mathord{\begin{tikzpicture}[baseline = -1mm]
  \draw[-,thick,darkblue] (0,0.2) to[out=180,in=90] (-.2,0);
  \draw[-,thick,darkblue] (0.2,0) to[out=90,in=0] (0,.2);
 \draw[-,thick,darkblue] (-.2,0) to[out=-90,in=180] (0,-0.2);
  \draw[-,thick,darkblue] (0,-0.2) to[out=0,in=-90] (0.2,0);
\end{tikzpicture}
}$ to the scalar $-z_1$, i.e., $\Psi_f$ factors through the specialization
$\H_k(\delta)$ where $\delta$ is the sum of the roots of the
polynomial $f(u)$.

The natural transformations $\Psi_f(c)$ and $\Psi_f(d)$ in the
previous paragraph come from the units and counits of the canonical
adjunctions making $(\ind^{n+1}_n,\res^{n+1}_n)$ into adjoint pairs.
In view of Theorem~\ref{thm2}(ii), we also get canonical adjunctions
the other
way around,
with units and counits defined by $\Psi_f(c')$ and $\Psi_f(d')$, respectively.
Thus, the induction and restriction functors $\ind_{n}^{n+1}$ and
$\res_n^{n+1}$ are biadjoint; see also \cite[Corollary 7.7.5]{Kbook} and
\cite[Proposition 5.13]{MS}.

One reason that cyclotomic quotients of the degenerate affine Hecke algebra are
important is that they can be used to
realize the 
{\em minimal categorifications} of
integrable lowest weight modules for the Lie algebra
$\mathfrak{g} := \mathfrak{sl}_\infty$ (if $\k$ is a field of characteristic $ 0$) or
$\mathfrak{g} := \widehat{\mathfrak{sl}}_p$ (if $\k$ is a field of characteristic $p > 0$),
e.g., see 
\cite{Ariki, BK}.
The following theorem shows that these minimal categorifications can be realized 
instead as {\em cyclotomic quotients of Heisenberg categories}.
All of this
should be compared with \cite[$\S$5.1.2]{Rou} (and \cite[Theorem 4.25]{Rou2}),
where the minimal categorification is
realized as a cyclotomic quotient of the corresponding Kac-Moody
2-category.
In the special case $\ell=1$, some closely related
constructions can be found in \cite{QSY}.
 
\begin{theorem}\label{darn}
Fix $f(u) = u^\ell + z_1 u^{\ell-1} + \cdots +
z_\ell \in \k[u]$ of degree $\ell=-k > 0$ as in
(\ref{wolf}). Let $\mathcal I_{f,1}$ be the $\k$-linear left tensor ideal of $\H_k$
generated by $f(x):\E\rightarrow \E$;
equivalently, by Lemma~\ref{belo} below, $\mathcal I_{f,1}$ is the
$\k$-linear left tensor
ideal generated by $1_\F:\F\rightarrow \F$ and $\mathord{\begin{tikzpicture}[baseline = -1mm]
  \draw[<-,thick,darkblue] (0,0.2) to[out=180,in=90] (-.2,0);
  \draw[-,thick,darkblue] (0.2,0) to[out=90,in=0] (0,.2);
 \draw[-,thick,darkblue] (-.2,0) to[out=-90,in=180] (0,-0.2);
  \draw[-,thick,darkblue] (0,-0.2) to[out=0,in=-90] (0.2,0);
   \node at (-0.68,0) {$\color{darkblue}\scriptstyle{r+k-1}$};
      \node at (-0.2,0) {$\dot$};
\end{tikzpicture}
}
 + z_r1_\unit:\unit\rightarrow\unit$
for $r=1,\dots,\ell$.
Let
$$
\ev:\mathcal{E}\!nd_\k\left(\bigoplus_{n \geq 0}  H_n^f\Mod\right) \rightarrow
\bigoplus_{n \geq 0} H^f_n\Mod$$ 
be the functor
defined by evaluating on the one-dimensional $H_0^f$-module.
Then $\ev \circ \Psi_f$
factors through the quotient category $\H_{f,1} := \H_k / \mathcal I_{f,1}$
to induce an equivalence of categories
$$
\psi_f:\Kar(\H_{f,1}) \rightarrow 
\bigoplus_{n \geq 0} H_n^f\proj,
$$
where $\Kar$ denotes additive Karoubi envelope and 
$\operatorname{pmod}$ denotes finitely generated projectives.
\end{theorem}

To get the full structure of a $\mathfrak{g}$-categorification
on $\Kar(\H_{f,1})$ in the sense of \cite[Definition 5.29]{Rou},
one also needs the endofunctors $E$ and $F$ defined by tensoring with
$\up$ and $\down$, respectively. 
Under the equivalence in Theorem~\ref{darn},
these
correspond to the induction and restriction functors on $\bigoplus_{n
  \geq 0} H_n^f\proj$. 
It is immediate from the definition of $\H_{f,1}$ 
that $E$ and $F$ are biadjoint and that the powers of $E$ 
admit the appropriate action of the degenerate
affine Hecke algebra. It just remains to check that 
the complexified Grothendieck group $\C\otimes_{\Z}
K_0(\Kar(\H_{f,1}))$ is the appropriate integrable representation of
$\g$. This follows from \cite{Ariki,BK}
using the equivalence in the theorem.

In \cite{Web}, Webster introduced 
{\em generalized cyclotomic quotients of Kac-Moody 2-categories} which
categorify lowest-tensored-highest 
weight representations; see also
\cite[$\S$4.2]{BD}. For $\mathfrak{sl}_\infty$ or
$\widehat{\mathfrak{sl}}_p$, Webster's categories can
also be realized as 
{\em generalized cyclotomic
quotients of Heisenberg categories}.
This will be explained elsewhere, but we can 
at least formulate the definition of these generalized cyclotomic quotients here.
Fix a pair of monic
polynomials\footnote{We stress here that $f'(u)$ denotes a different polynomial; it is {\em not} the derivative of $f(u)$!}
$f(u), f'(u) \in \k[u]$ of degrees $\ell, \ell'\geq 0$, respectively,
and define $k := \ell'-\ell$ and $\delta_r, \delta_r' \in \k$ so that
\begin{align}\label{delta}
\delta(u) &= \delta_0 + \delta_1 u^{-1} + \delta_2 u^{-2} + \cdots := u^{-k} f'(u) / f(u)
\in \k[[u^{-1}]],\\
\delta'(u) &= \delta'_0 + \delta'_1 u^{-1} + \delta'_2 u^{-2} + \cdots := -u^k f(u) / f'(u)
\in \k[[u^{-1}]].\label{deltap}
\end{align}
Then the corresponding generalized cyclotomic quotient of $\H_k$ is the
$\k$-linear category
\begin{equation}\label{blake}
\H_{f,f'} := \H_k /
\mathcal I_{f,f'}
\end{equation}
where $\mathcal I_{f,f'}$ is the $\k$-linear left tensor ideal of
$\H_k$
generated by $f(x):\E\rightarrow \E$ and $\mathord{\begin{tikzpicture}[baseline = -1mm]
  \draw[-,thick,darkblue] (0,0.2) to[out=180,in=90] (-.2,0);
  \draw[->,thick,darkblue] (0.2,0) to[out=90,in=0] (0,.2);
 \draw[-,thick,darkblue] (-.2,0) to[out=-90,in=180] (0,-0.2);
  \draw[-,thick,darkblue] (0,-0.2) to[out=0,in=-90] (0.2,0);
   \node at (0.68,0) {$\color{darkblue}\scriptstyle{r-k-1}$};
      \node at (0.2,0) {$\dot$};
\end{tikzpicture}
}-\delta_r 1_\unit:\unit\rightarrow\unit$ for $r = 1,\dots,\ell'$.
These categories were introduced already in the case that $\ell=\ell'$
in \cite{BCNR}.

\begin{lemma}\label{belo}
The ideal $\mathcal I_{f,f'}$ 
can be defined equivalently as the $\k$-linear left tensor ideal of $\H_k$ generated by
$f'(x'):\F\rightarrow \F$ and $
\mathord{\begin{tikzpicture}[baseline = -1mm]
  \draw[<-,thick,darkblue] (0,0.2) to[out=180,in=90] (-.2,0);
  \draw[-,thick,darkblue] (0.2,0) to[out=90,in=0] (0,.2);
 \draw[-,thick,darkblue] (-.2,0) to[out=-90,in=180] (0,-0.2);
  \draw[-,thick,darkblue] (0,-0.2) to[out=0,in=-90] (0.2,0);
   \node at (-0.68,0) {$\color{darkblue}\scriptstyle{r+k-1}$};
      \node at (-0.2,0) {$\dot$};
\end{tikzpicture}
}
-\delta_r' 1_\unit:\unit\rightarrow\unit$ for $r = 1,\dots,\ell$.
It also
contains 
$\mathord{\begin{tikzpicture}[baseline = -1mm]
  \draw[-,thick,darkblue] (0,0.2) to[out=180,in=90] (-.2,0);
  \draw[->,thick,darkblue] (0.2,0) to[out=90,in=0] (0,.2);
 \draw[-,thick,darkblue] (-.2,0) to[out=-90,in=180] (0,-0.2);
  \draw[-,thick,darkblue] (0,-0.2) to[out=0,in=-90] (0.2,0);
   \node at (0.68,0) {$\color{darkblue}\scriptstyle{r-k-1}$};
      \node at (0.2,0) {$\dot$};
\end{tikzpicture}
}-\delta_r 1_\unit$
and
$
\mathord{\begin{tikzpicture}[baseline = -1mm]
  \draw[<-,thick,darkblue] (0,0.2) to[out=180,in=90] (-.2,0);
  \draw[-,thick,darkblue] (0.2,0) to[out=90,in=0] (0,.2);
 \draw[-,thick,darkblue] (-.2,0) to[out=-90,in=180] (0,-0.2);
  \draw[-,thick,darkblue] (0,-0.2) to[out=0,in=-90] (0.2,0);
   \node at (-0.68,0) {$\color{darkblue}\scriptstyle{r+k-1}$};
      \node at (-0.2,0) {$\dot$};
\end{tikzpicture}
}
-\delta_r' 1_\unit$ for all $r \geq 0$.
\end{lemma}

Let us finally mention that there is also a quantum analog $\H_k(z,t)$
of the
Heisenberg category $\H_k(\delta)$. This will be defined in a sequel to this
article \cite{B2}.
Even in the case that $k=-1$, our approach is different
to that of \cite{LS} as we require that the polynomial generator $X$
is invertible, i.e., we 
incorporate the entire affine Hecke algebra into
the definition (rather than the $q$-deformed degenerate affine Hecke
algebra used in \cite{LS}).
The quantum Heisenberg category $\H_0(z,t)$ of
central charge zero is the {\em affine oriented skein category}
$\mathcal{AOS}(z,t)$ from
\cite[$\S$4]{Bnew}.
Further generalizations incorporating Clifford 
and Frobenius
superalgebras into the definition have also recently emerged building
on the approach taken in this article; see
\cite{CK}
(which extends \cite{HS} to arbitrary central charge) 
and \cite{S} (which extends \cite{RS}).

\vspace{2mm}
\noindent
{\em Acknowledgements.}
My thanks go to Jonathan Comes for 
correcting some sign errors in the
first version of the article.

\section{Analysis of the inversion relation}

This section is the technical heart of the paper.
The development is similar
to that of \cite{B} but with subtlely different signs.
Going back to the original definition of $\H_k$ from
Definition~\ref{maindef}, we begin our study by {\em defining} the
 {downward dots and crossings} to be
the right mates of the upward dots and crossings:
\begin{align}
&
x' =\mathord{
\begin{tikzpicture}[baseline = 0]
	\draw[<-,thick,darkblue] (0.08,-.3) to (0.08,.4);
      \node at (0.08,0.1) {$\dot$};
\end{tikzpicture}
}
:=
\mathord{
\begin{tikzpicture}[baseline = 0]
  \draw[->,thick,darkblue] (0.3,0) to (0.3,-.4);
	\draw[-,thick,darkblue] (0.3,0) to[out=90, in=0] (0.1,0.4);
	\draw[-,thick,darkblue] (0.1,0.4) to[out = 180, in = 90] (-0.1,0);
	\draw[-,thick,darkblue] (-0.1,0) to[out=-90, in=0] (-0.3,-0.4);
	\draw[-,thick,darkblue] (-0.3,-0.4) to[out = 180, in =-90] (-0.5,0);
  \draw[-,thick,darkblue] (-0.5,0) to (-0.5,.4);
   \node at (-0.1,0) {$\dot$};
\end{tikzpicture}
},
\hspace{20mm}
s' = \mathord{
\begin{tikzpicture}[baseline = 0]
	\draw[<-,thick,darkblue] (0.28,-.3) to (-0.28,.4);
	\draw[<-,thick,darkblue] (-0.28,-.3) to (0.28,.4);
\end{tikzpicture}
}:=
\mathord{
\begin{tikzpicture}[baseline = 0]
	\draw[<-,thick,darkblue] (0.3,-.5) to (-0.3,.5);
	\draw[-,thick,darkblue] (-0.2,-.2) to (0.2,.3);
        \draw[-,thick,darkblue] (0.2,.3) to[out=50,in=180] (0.5,.5);
        \draw[->,thick,darkblue] (0.5,.5) to[out=0,in=90] (0.8,-.5);
        \draw[-,thick,darkblue] (-0.2,-.2) to[out=230,in=0] (-0.5,-.5);
        \draw[-,thick,darkblue] (-0.5,-.5) to[out=180,in=-90] (-0.8,.5);
\end{tikzpicture}
}\:.
\end{align}
The following relations are immediate from these definitions:
\begin{align}\label{cough1}
\mathord{
\begin{tikzpicture}[baseline = 0]
	\draw[<-,thick,darkblue] (0.4,0.4) to[out=-90, in=0] (0.1,-0.1);
	\draw[-,thick,darkblue] (0.1,-0.1) to[out = 180, in = -90] (-0.2,0.4);
      \node at (0.37,0.15) {$\dot$};
\end{tikzpicture}
}
&=
\mathord{
\begin{tikzpicture}[baseline = 0]
	\draw[<-,thick,darkblue] (0.4,0.4) to[out=-90, in=0] (0.1,-0.1);
	\draw[-,thick,darkblue] (0.1,-0.1) to[out = 180, in = -90] (-0.2,0.4);
      \node at (-0.16,0.15) {$\dot$};
\end{tikzpicture}
},&
\mathord{
\begin{tikzpicture}[baseline = 0]
\draw[->,thick,darkblue](-.5,-.3) to (0,.4);
	\draw[<-,thick,darkblue] (0.3,0.4) to[out=-90, in=0] (0,-0.1);
	\draw[-,thick,darkblue] (0,-0.1) to[out = 180, in = -40] (-0.5,0.4);
\end{tikzpicture}
}&=
\mathord{
\begin{tikzpicture}[baseline = 0]
\draw[->,thick,darkblue](.6,-.3) to (.1,.4);
	\draw[<-,thick,darkblue] (0.6,0.4) to[out=-140, in=0] (0.1,-0.1);
	\draw[-,thick,darkblue] (0.1,-0.1) to[out = 180, in = -90] (-0.2,0.4);
\end{tikzpicture}
},&
\mathord{
\begin{tikzpicture}[baseline = 0]
\draw[<-,thick,darkblue](.6,-.3) to (.1,.4);
	\draw[<-,thick,darkblue] (0.6,0.4) to[out=-140, in=0] (0.1,-0.1);
	\draw[-,thick,darkblue] (0.1,-0.1) to[out = 180, in = -90] (-0.2,0.4);
\end{tikzpicture}
}&=
\mathord{
\begin{tikzpicture}[baseline = 0]
\draw[<-,thick,darkblue](-.5,-.3) to (0,.4);
	\draw[<-,thick,darkblue] (0.3,0.4) to[out=-90, in=0] (0,-0.1);
	\draw[-,thick,darkblue] (0,-0.1) to[out = 180, in = -40] (-0.5,0.4);
\end{tikzpicture}
},
\\\label{cough2}
\mathord{
\begin{tikzpicture}[baseline = 1mm]
	\draw[<-,thick,darkblue] (0.4,-0.1) to[out=90, in=0] (0.1,0.4);
	\draw[-,thick,darkblue] (0.1,0.4) to[out = 180, in = 90] (-0.2,-0.1);
      \node at (-0.14,0.25) {$\dot$};
\end{tikzpicture}
}&
=
\mathord{
\begin{tikzpicture}[baseline = 1mm]
	\draw[<-,thick,darkblue] (0.4,-0.1) to[out=90, in=0] (0.1,0.4);
	\draw[-,thick,darkblue] (0.1,0.4) to[out = 180, in = 90] (-0.2,-0.1);
      \node at (0.32,0.25) {$\dot$};
\end{tikzpicture}
},&
\mathord{
\begin{tikzpicture}[baseline = 0]
\draw[<-,thick,darkblue](.6,.4) to (.1,-.3);
	\draw[<-,thick,darkblue] (0.6,-0.3) to[out=140, in=0] (0.1,0.2);
	\draw[-,thick,darkblue] (0.1,0.2) to[out = -180, in = 90] (-0.2,-0.3);
\end{tikzpicture}
}&=
\mathord{
\begin{tikzpicture}[baseline = 0]
\draw[<-,thick,darkblue](-.5,.4) to (0,-.3);
	\draw[<-,thick,darkblue] (0.3,-0.3) to[out=90, in=0] (0,0.2);
	\draw[-,thick,darkblue] (0,0.2) to[out = -180, in = 40] (-0.5,-0.3);
\end{tikzpicture}
},&
\mathord{
\begin{tikzpicture}[baseline = 0]
\draw[->,thick,darkblue](-.5,.4) to (0,-.3);
	\draw[<-,thick,darkblue] (0.3,-0.3) to[out=90, in=0] (0,0.2);
	\draw[-,thick,darkblue] (0,0.2) to[out = -180, in = 40] (-0.5,-0.3);
\end{tikzpicture}
}&=
\mathord{
\begin{tikzpicture}[baseline = 0]
\draw[->,thick,darkblue](.6,.4) to (.1,-.3);
	\draw[<-,thick,darkblue] (0.6,-0.3) to[out=140, in=0] (0.1,0.2);
	\draw[-,thick,darkblue] (0.1,0.2) to[out = -180, in = 90] (-0.2,-0.3);
\end{tikzpicture}
}.
\end{align}
Also, the following relations are easily deduced by attaching
rightward cups and caps to the degenerate Hecke relations, then ``rotating'' the pictures using the
definitions of the rightwards/downwards crossings and the downwards 
dots:
\begin{align}\label{spit1}
\mathord{
\begin{tikzpicture}[baseline = 0]
	\draw[-,thick,darkblue] (0.28,0) to[out=90,in=-90] (-0.28,.6);
	\draw[-,thick,darkblue] (-0.28,0) to[out=90,in=-90] (0.28,.6);
	\draw[<-,thick,darkblue] (0.28,-.6) to[out=90,in=-90] (-0.28,0);
	\draw[<-,thick,darkblue] (-0.28,-.6) to[out=90,in=-90] (0.28,0);
\end{tikzpicture}
}
&=\mathord{
\begin{tikzpicture}[baseline = 0]
	\draw[<-,thick,darkblue] (0.08,-.6) to (0.08,.6);
	\draw[<-,thick,darkblue] (-0.28,-.6) to (-0.28,.6);
 \end{tikzpicture}
},
&\mathord{
\begin{tikzpicture}[baseline = 0]
	\draw[<-,thick,darkblue] (0.45,.6) to (-0.45,-.6);
	\draw[<-,thick,darkblue] (0.45,-.6) to (-0.45,.6);
        \draw[-,thick,darkblue] (0,-.6) to[out=90,in=-90] (.45,0);
        \draw[->,thick,darkblue] (0.45,0) to[out=90,in=-90] (0,0.6);
\end{tikzpicture}
}&=
\mathord{
\begin{tikzpicture}[baseline = 0]
	\draw[<-,thick,darkblue] (0.45,.6) to (-0.45,-.6);
	\draw[<-,thick,darkblue] (0.45,-.6) to (-0.45,.6);
        \draw[-,thick,darkblue] (0,-.6) to[out=90,in=-90] (-.45,0);
        \draw[->,thick,darkblue] (-0.45,0) to[out=90,in=-90] (0,0.6);
\end{tikzpicture}
}\:,
&
\mathord{
\begin{tikzpicture}[baseline = 0]
	\draw[->,thick,darkblue] (0.45,.6) to (-0.45,-.6);
	\draw[<-,thick,darkblue] (0.45,-.6) to (-0.45,.6);
        \draw[<-,thick,darkblue] (0,-.6) to[out=90,in=-90] (.45,0);
        \draw[-,thick,darkblue] (0.45,0) to[out=90,in=-90] (0,0.6);
\end{tikzpicture}
}
&=\mathord{
\begin{tikzpicture}[baseline = 0]
	\draw[->,thick,darkblue] (0.45,.6) to (-0.45,-.6);
	\draw[<-,thick,darkblue] (0.45,-.6) to (-0.45,.6);
        \draw[<-,thick,darkblue] (0,-.6) to[out=90,in=-90] (-.45,0);
        \draw[-,thick,darkblue] (-0.45,0) to[out=90,in=-90] (0,0.6);
\end{tikzpicture}
},\end{align}
\begin{align}
\label{spit2}
\mathord{
\begin{tikzpicture}[baseline = 1mm]
	\draw[<-,thick,darkblue] (0.3,.6) to (-0.3,-.2);
	\draw[<-,thick,darkblue] (0.3,-.2) to (-0.3,.6);
     \node at (0.16,0.42) {$\dot$};
 \end{tikzpicture}
}
-
\mathord{
\begin{tikzpicture}[baseline = 1mm]
	\draw[<-,thick,darkblue] (0.3,.6) to (-0.3,-.2);
	\draw[<-,thick,darkblue] (0.3,-.2) to (-0.3,.6);
      \node at (-0.16,-0.02) {$\dot$};
\end{tikzpicture}
}
&=
\mathord{
\begin{tikzpicture}[baseline = 1mm]
 	\draw[-,thick,darkblue] (-0.33,-.2) to [out=90,in=180] (-0.1,0.1);
 	\draw[->,thick,darkblue] (-0.1,0.1) to [out=0,in=90] (0.13,-.2);
	\draw[-,thick,darkblue] (-0.33,.6) to [out=-90,in=180] (-0.1,.3);
	\draw[->,thick,darkblue] (-0.1,.3) to [out=0,in=-90](0.13,.6);
\end{tikzpicture}
}
=
\mathord{\begin{tikzpicture}[baseline = 1mm]
 	\draw[<-,thick,darkblue] (0.3,.6) to (-0.3,-.2);
	\draw[<-,thick,darkblue] (0.3,-.2) to (-0.3,.6);
      \node at (-0.16,0.42) {$\dot$};
\end{tikzpicture}
}
-\mathord{
\begin{tikzpicture}[baseline = 1mm]
	\draw[<-,thick,darkblue] (0.3,.6) to (-0.3,-.2);
	\draw[<-,thick,darkblue] (0.3,-.2) to (-0.3,.6);
      \node at (0.16,-0.02) {$\dot$};
\end{tikzpicture}
}
,
&\mathord{\begin{tikzpicture}[baseline = 1mm]
	\draw[->,thick,darkblue] (0.3,.6) to (-0.3,-.2);
	\draw[<-,thick,darkblue] (0.3,-.2) to (-0.3,.6);
      \node at (0.16,.42) {$\dot$};
\end{tikzpicture}
}
-\mathord{
\begin{tikzpicture}[baseline = 1mm]
	\draw[->,thick,darkblue] (0.3,.6) to (-0.3,-.2);
	\draw[<-,thick,darkblue] (0.3,-.2) to (-0.3,.6);
      \node at (-0.16,0) {$\dot$};
\end{tikzpicture}
}
&=
\mathord{
\begin{tikzpicture}[baseline = 1mm]
 	\draw[<-,thick,darkblue] (0.08,-.2) to (0.08,.6);
	\draw[<-,thick,darkblue] (-0.28,-.2) to (-0.28,.6);
\end{tikzpicture}
}
=
\mathord{\begin{tikzpicture}[baseline =1mm]
	\draw[->,thick,darkblue] (0.3,.6) to (-0.3,-.2);
	\draw[<-,thick,darkblue] (0.3,-.2) to (-0.3,.6);
      \node at (0.16,0) {$\dot$};
\end{tikzpicture}
}
-\mathord{
\begin{tikzpicture}[baseline = 1mm]
	\draw[->,thick,darkblue] (0.3,.6) to (-0.3,-.2);
	\draw[<-,thick,darkblue] (0.3,-.2) to (-0.3,.6);
      \node at (-0.16,0.42) {$\dot$};
\end{tikzpicture}
}.
\end{align}

The important symmetry $\omega$ constructed in the next lemma is
often useful since it reduces to the case that $k \geq 0$.
In words, $\omega$ reflects in a horizontal axis then
multiplies by $(-1)^{x}$, where $x$ is the total number of crossings
appearing in the diagram.
This heuristic also holds for all of the other morphisms
defined diagrammatically 
below, but in general the sign becomes $(-1)^{x+y}$ where $x$ is the total
number of crossings and $y$ is the total number of 
leftward
cups and caps (not counting the decorated caps and cups to be
introduced shortly which are
labelled with the symbol $\color{darkblue}\diamond$).

\begin{lemma}\label{omega}
There is an isomorphism of monoidal categories
$\omega:\H_k \stackrel{\sim}{\rightarrow} \H_{-k}^{\operatorname{op}}$
switching the objects $\E$ and $\F$, and defined on generating morphisms
by $x \mapsto x', s \mapsto -s', c \mapsto d$ and $d \mapsto c$.
\end{lemma}

\proof
The existence of $\omega$ follows by
a straightforward relation check. Use (\ref{spit1})--(\ref{spit2}) for
the degenerate Hecke relations.
The need to switch $k$ and $-k$ comes from the inversion relations.
To see that $\omega$ is an isomorphism, notice by the 
right adjunction 
relations that $\omega(x') = x$ and $\omega(s') = -s$, hence, 
$\omega^2 = \id$.
\endproof

The inversion relation means that
there are some
as yet unnamed
generating morphisms in $\H_k$ which are
the matrix entries of
two-sided inverses to the morphism (\ref{invrel1})--(\ref{invrel2}).
We next introduce notation for these matrix entries.
First define
\begin{align*}
t'=\mathord{
\begin{tikzpicture}[baseline = 0]
	\draw[->,thick,darkblue] (0.28,-.3) to (-0.28,.4);
	\draw[<-,thick,darkblue] (-0.28,-.3) to (0.28,.4);
\end{tikzpicture}
}
:\F\otimes \E \rightarrow \E\otimes \F,
\end{align*}
and the decorated leftward cups and caps
\begin{align*}
\mathord{
\begin{tikzpicture}[baseline = 1mm]
	\draw[-,thick,darkblue] (0.4,0.4) to[out=-90, in=0] (0.1,0);
	\draw[->,thick,darkblue] (0.1,0) to[out = 180, in = -90] (-0.2,0.4);
      \node at (0.12,-0.2) {$\color{darkblue}\scriptstyle{r}$};
      \node at (0.12,0.01) {$\diamond$};
\end{tikzpicture}
}
&:\unit \rightarrow \E\otimes \F,
&
\mathord{
\begin{tikzpicture}[baseline = 1mm]
	\draw[-,thick,darkblue] (0.4,0) to[out=90, in=0] (0.1,0.4);
	\draw[->,thick,darkblue] (0.1,0.4) to[out = 180, in = 90] (-0.2,0);
      \node at (0.12,0.6) {$\color{darkblue}\scriptstyle{r}$};
      \node at (0.12,0.4) {$\diamond$};
\end{tikzpicture}
}
&:\F\otimes \E \rightarrow \unit
\end{align*}
for $0 \leq r < k$ or $0 \leq r < -k$, respectively, by declaring that 
\begin{align}\label{cutting1}
\left[
\displaystyle\mathord{
\begin{tikzpicture}[baseline = 0]
	\draw[->,thick,darkblue] (0.28,-.3) to (-0.28,.4);
	\draw[<-,thick,darkblue] (-0.28,-.3) to (0.28,.4);
\end{tikzpicture}
}\:\:\:
\mathord{
\begin{tikzpicture}[baseline = 1mm]
	\draw[-,thick,darkblue] (0.4,0.4) to[out=-90, in=0] (0.1,0);
	\draw[->,thick,darkblue] (0.1,0) to[out = 180, in = -90] (-0.2,0.4);
      \node at (0.12,-0.2) {$\color{darkblue}\scriptstyle{0}$};
      \node at (0.12,0.01) {$\diamond$};
\end{tikzpicture}
}
\:\:\:
\mathord{
\begin{tikzpicture}[baseline = 1mm]
	\draw[-,thick,darkblue] (0.4,0.4) to[out=-90, in=0] (0.1,0);
	\draw[->,thick,darkblue] (0.1,0) to[out = 180, in = -90] (-0.2,0.4);
      \node at (0.12,-0.2) {$\color{darkblue}\scriptstyle{1}$};
      \node at (0.12,0.01) {$\diamond$};
\end{tikzpicture}
}
\:\:\:\cdots\:\:\:
\mathord{
\begin{tikzpicture}[baseline = 1mm]
	\draw[-,thick,darkblue] (0.4,0.4) to[out=-90, in=0] (0.1,0);
	\draw[->,thick,darkblue] (0.1,0) to[out = 180, in = -90] (-0.2,0.4);
      \node at (0.12,-0.2) {$\color{darkblue}\scriptstyle{k-1}$};
      \node at (0.12,0.01) {$\diamond$};
\end{tikzpicture}
}
\right]
&:=
\left(\left[
\mathord{
\begin{tikzpicture}[baseline = 0]
	\draw[<-,thick,darkblue] (0.28,-.3) to (-0.28,.4);
	\draw[->,thick,darkblue] (-0.28,-.3) to (0.28,.4);
\end{tikzpicture}
}\:\:\:
\mathord{
\begin{tikzpicture}[baseline = 1mm]
	\draw[<-,thick,darkblue] (0.4,0) to[out=90, in=0] (0.1,0.4);
	\draw[-,thick,darkblue] (0.1,0.4) to[out = 180, in = 90] (-0.2,0);
\end{tikzpicture}
}
\:\:\:
\mathord{
\begin{tikzpicture}[baseline = 1mm]
	\draw[<-,thick,darkblue] (0.4,0) to[out=90, in=0] (0.1,0.4);
	\draw[-,thick,darkblue] (0.1,0.4) to[out = 180, in = 90] (-0.2,0);
      \node at (-0.15,0.2) {$\dot$};
\end{tikzpicture}
}
\:\:\:\cdots
\:\:\:\mathord{
\begin{tikzpicture}[baseline = 1mm]
	\draw[<-,thick,darkblue] (0.4,0) to[out=90, in=0] (0.1,0.4);
	\draw[-,thick,darkblue] (0.1,0.4) to[out = 180, in = 90] (-0.2,0);
      \node at (-0.53,0.2) {$\color{darkblue}\scriptstyle{k-1}$};
      \node at (-0.15,0.2) {$\dot$};
\end{tikzpicture}
}\right]^{T}\right)^{-1}\\\intertext{if $k \geq 0$, or}
\left[
\mathord{
\begin{tikzpicture}[baseline = 0]
	\draw[->,thick,darkblue] (0.28,-.3) to (-0.28,.4);
	\draw[<-,thick,darkblue] (-0.28,-.3) to (0.28,.4);
\end{tikzpicture}
}\:\:\:
\mathord{
\begin{tikzpicture}[baseline = 1mm]
	\draw[-,thick,darkblue] (0.4,0) to[out=90, in=0] (0.1,0.4);
	\draw[->,thick,darkblue] (0.1,0.4) to[out = 180, in = 90] (-0.2,0);
      \node at (0.12,0.6) {$\color{darkblue}\scriptstyle{0}$};
      \node at (0.12,0.4) {$\diamond$};
\end{tikzpicture}
}
\:\:\:
\mathord{
\begin{tikzpicture}[baseline = 1mm]
	\draw[-,thick,darkblue] (0.4,0) to[out=90, in=0] (0.1,0.4);
	\draw[->,thick,darkblue] (0.1,0.4) to[out = 180, in = 90] (-0.2,0);
      \node at (0.12,0.6) {$\color{darkblue}\scriptstyle{1}$};
      \node at (0.12,0.4) {$\diamond$};
\end{tikzpicture}
}
\:\:\:\cdots\:\:\:
\mathord{
\begin{tikzpicture}[baseline = 1mm]
	\draw[-,thick,darkblue] (0.4,0) to[out=90, in=0] (0.1,0.4);
	\draw[->,thick,darkblue] (0.1,0.4) to[out = 180, in = 90] (-0.2,0);
      \node at (0.12,0.6) {$\color{darkblue}\scriptstyle{-k-1}$};
      \node at (0.12,0.4) {$\diamond$};
\end{tikzpicture}
}
\right]^T\!\!\!\!&:=
\left[
\mathord{
\begin{tikzpicture}[baseline = 0]
	\draw[<-,thick,darkblue] (0.28,-.3) to (-0.28,.4);
	\draw[->,thick,darkblue] (-0.28,-.3) to (0.28,.4);
\end{tikzpicture}
}\:\:\:
\mathord{
\begin{tikzpicture}[baseline =1mm]
	\draw[<-,thick,darkblue] (0.4,0.4) to[out=-90, in=0] (0.1,0);
	\draw[-,thick,darkblue] (0.1,0) to[out = 180, in = -90] (-0.2,0.4);
\end{tikzpicture}
}
\:\:\:\mathord{
\begin{tikzpicture}[baseline =1mm]
	\draw[<-,thick,darkblue] (0.4,0.4) to[out=-90, in=0] (0.1,0);
	\draw[-,thick,darkblue] (0.1,0) to[out = 180, in = -90] (-0.2,0.4);
      \node at (0.37,0.2) {$\dot$};
\end{tikzpicture}
}
\:\:\:\cdots\:\:\:\mathord{
\begin{tikzpicture}[baseline =1mm]
	\draw[<-,thick,darkblue] (0.4,0.4) to[out=-90, in=0] (0.1,0);
	\draw[-,thick,darkblue] (0.1,0) to[out = 180, in = -90] (-0.2,0.4);
      \node at (0.77,0.2) {$\color{darkblue}\scriptstyle{-k-1}$};
      \node at (0.37,0.2) {$\dot$};
\end{tikzpicture}
}
\right]^{-1}\label{cutting2}
\end{align}
if $k < 0$.
Then 
we set
\begin{align}\label{c}
c'&=\mathord{
\begin{tikzpicture}[baseline = 1mm]
	\draw[-,thick,darkblue] (0.4,0.4) to[out=-90, in=0] (0.1,0);
	\draw[->,thick,darkblue] (0.1,0) to[out = 180, in = -90] (-0.2,0.4);
 \end{tikzpicture}
}\!\!\!\!
&:=
\left\{
\begin{array}{rl}
\!\!\! -\:\mathord{
\begin{tikzpicture}[baseline = 0.5mm]
	\draw[-,thick,darkblue] (0.4,0.4) to[out=-90, in=0] (0.1,0);
	\draw[->,thick,darkblue] (0.1,0) to[out = 180, in = -90] (-0.2,0.4);
      \node at (0.16,-0.21) {$\color{darkblue}\scriptstyle{k-1}$};
      \node at (0.12,0.01) {$\diamond$};
\end{tikzpicture}
}
\hspace{3.5mm}&\!\!\!\text{if $k > 0$,}\\
\mathord{
\begin{tikzpicture}[baseline = 0]
	\draw[-,thick,darkblue] (0.28,.6) to[out=240,in=90] (-0.28,-.1);
	\draw[<-,thick,darkblue] (-0.28,.6) to[out=300,in=90] (0.28,-0.1);
	\draw[-,thick,darkblue] (0.28,-0.1) to[out=-90, in=0] (0,-0.4);
	\draw[-,thick,darkblue] (0,-0.4) to[out = 180, in = -90] (-0.28,-0.1);
      \node at (0.54,-0.1) {$\color{darkblue}\scriptstyle{-k}$};
      \node at (0.27,-0.1) {$\dot$};
\end{tikzpicture}
}
&\!\!\!\text{if $k \leq 0$,}
\end{array}\right.
&
\:\:\:\:\:\:d'=\mathord{
\begin{tikzpicture}[baseline = 1mm]
	\draw[-,thick,darkblue] (0.4,0) to[out=90, in=0] (0.1,0.4);
	\draw[->,thick,darkblue] (0.1,0.4) to[out = 180, in = 90] (-0.2,0);
\end{tikzpicture}
}
&\!\!\!\!:=
\left\{
\begin{array}{rl}
\!\!\mathord{
\begin{tikzpicture}[baseline = -1.5mm]
	\draw[-,thick,darkblue] (0.28,-.6) to[out=120,in=-90] (-0.28,.1);
	\draw[<-,thick,darkblue] (-0.28,-.6) to[out=60,in=-90] (0.28,.1);
	\draw[-,thick,darkblue] (0.28,0.1) to[out=90, in=0] (0,0.4);
	\draw[-,thick,darkblue] (0,0.4) to[out = 180, in = 90] (-0.28,0.1);
      \node at (-0.44,0.15) {$\color{darkblue}\scriptstyle{k}$};
      \node at (-0.27,0.15) {$\dot$};
\end{tikzpicture}
}
\hspace{4.5mm}
&\!\!\!\text{if $k  \geq 0$,}\\
\!\!\mathord{
\begin{tikzpicture}[baseline = 1mm]
	\draw[-,thick,darkblue] (0.4,0) to[out=90, in=0] (0.1,0.4);
	\draw[->,thick,darkblue] (0.1,0.4) to[out = 180, in = 90] (-0.2,0);
      \node at (0.06,0.62) {$\color{darkblue}\scriptstyle{-k-1}$};
      \node at (0.12,0.4) {$\diamond$};
\end{tikzpicture}}
\hspace{3mm}&\!\!\!\text{if $k < 0$.}
\end{array}\right.
\end{align}
From these definitions, it follows that
\begin{align}\label{posneg}
\mathord{
\begin{tikzpicture}[baseline = 0]
	\draw[<-,thick,darkblue] (0.08,-.6) to (0.08,.6);
	\draw[->,thick,darkblue] (-0.28,-.6) to (-0.28,.6);
\end{tikzpicture}
}
=
\mathord{
\begin{tikzpicture}[baseline = 0mm]
	\draw[->,thick,darkblue] (0.28,0) to[out=90,in=-90] (-0.28,.6);
	\draw[-,thick,darkblue] (-0.28,0) to[out=90,in=-90] (0.28,.6);
	\draw[<-,thick,darkblue] (0.28,-.6) to[out=90,in=-90] (-0.28,0);
	\draw[-,thick,darkblue] (-0.28,-.6) to[out=90,in=-90] (0.28,0);
\end{tikzpicture}
}
+\sum_{r=0}^{k-1}
\mathord{
\begin{tikzpicture}[baseline=0mm]
	\draw[-,thick,darkblue] (0.3,0.6) to[out=-90, in=0] (0,0.1);
	\draw[->,thick,darkblue] (0,0.1) to[out = 180, in = -90] (-0.3,0.6);
      \node at (0.02,0.28) {$\color{darkblue}\scriptstyle{r}$};
	\draw[<-,thick,darkblue] (0.3,-.6) to[out=90, in=0] (0,-0.1);
	\draw[-,thick,darkblue] (0,-0.1) to[out = 180, in = 90] (-0.3,-.6);
   \node at (-0.25,-0.4) {$\dot$};
   \node at (-.4,-.3) {$\color{darkblue}\scriptstyle{r}$};
      \node at (0.02,0.11) {$\diamond$};
\end{tikzpicture}}
,\qquad\qquad
\mathord{
\begin{tikzpicture}[baseline = 0]
	\draw[->,thick,darkblue] (0.08,-.6) to (0.08,.6);
	\draw[<-,thick,darkblue] (-0.28,-.6) to (-0.28,.6);
\end{tikzpicture}
}=
\mathord{
\begin{tikzpicture}[baseline = 0mm]
	\draw[-,thick,darkblue] (0.28,0) to[out=90,in=-90] (-0.28,.6);
	\draw[->,thick,darkblue] (-0.28,0) to[out=90,in=-90] (0.28,.6);
	\draw[-,thick,darkblue] (0.28,-.6) to[out=90,in=-90] (-0.28,0);
	\draw[<-,thick,darkblue] (-0.28,-.6) to[out=90,in=-90] (0.28,0);
\end{tikzpicture}
}
+\sum_{r=0}^{-k-1}
\mathord{
\begin{tikzpicture}[baseline=0mm]
	\draw[<-,thick,darkblue] (0.3,0.6) to[out=-90, in=0] (0,.1);
	\draw[-,thick,darkblue] (0,.1) to[out = 180, in = -90] (-0.3,0.6);
      \node at (0.02,-0.28) {$\color{darkblue}\scriptstyle{r}$};
	\draw[-,thick,darkblue] (0.3,-.6) to[out=90, in=0] (0,-0.1);
	\draw[->,thick,darkblue] (0,-0.1) to[out = 180, in = 90] (-0.3,-.6);
   \node at (0.27,0.3) {$\dot$};
   \node at (.43,.3) {$\color{darkblue}\scriptstyle{r}$};
      \node at (0.02,-0.1) {$\diamond$};
\end{tikzpicture}},
\end{align}
with the right hand sides being sums of mutually orthogonal
idempotents.
Also
\begin{align}\label{life1}
\mathord{
\begin{tikzpicture}[baseline = -2mm]
	\draw[<-,thick,darkblue] (0.28,.4) to[out=240,in=90] (-0.28,-.3);
	\draw[-,thick,darkblue] (-0.28,.4) to[out=300,in=90] (0.28,-0.3);
  	\draw[-,thick,darkblue] (0.28,-0.3) to[out=-90, in=0] (0,-0.6);
	\draw[-,thick,darkblue] (0,-0.6) to[out = 180, in = -90] (-0.28,-0.3);
\end{tikzpicture}
}&=
\mathord{
\begin{tikzpicture}[baseline = -2mm]
	\draw[-,thick,darkblue] (0.28,-.6) to[out=120,in=-90] (-0.28,.1);
	\draw[<-,thick,darkblue] (-0.28,-.6) to[out=60,in=-90] (0.28,.1);
	\draw[-,thick,darkblue] (0.28,0.1) to[out=90, in=0] (0,0.4);
	\draw[-,thick,darkblue] (0,0.4) to[out = 180, in = 90] (-0.28,0.1);
      \node at (-0.45,0.05) {$\color{darkblue}\scriptstyle{r}$};
      \node at (-0.27,0.05) {$\dot$};
\end{tikzpicture}
}
=0
&
\text{and}
\qquad\qquad
\mathord{
\begin{tikzpicture}[baseline = 2mm]
  \draw[<-,thick,darkblue] (0,0.6) to[out=180,in=90] (-.3,0.3);
  \draw[-,thick,darkblue] (0.3,0.3) to[out=90,in=0] (0,.6);
 \draw[-,thick,darkblue] (-.3,0.3) to[out=-90,in=180] (0,0);
  \draw[-,thick,darkblue] (0,0) to[out=0,in=-90] (0.3,0.3);
   \node at (-0.3,0.3) {$\dot$};
   \node at (-0.5,0.3) {$\color{darkblue}\scriptstyle{r}$};
\end{tikzpicture}
}&= -\delta_{r,k-1} 1_\unit\\
\intertext{if $0 \leq r < k$, or}
\mathord{
\begin{tikzpicture}[baseline = -1mm]
	\draw[<-,thick,darkblue] (0.28,-.6) to[out=120,in=-90] (-0.28,.1);
	\draw[-,thick,darkblue] (-0.28,-.6) to[out=60,in=-90] (0.28,.1);
	\draw[-,thick,darkblue] (0.28,0.1) to[out=90, in=0] (0,0.4);
	\draw[-,thick,darkblue] (0,0.4) to[out = 180, in = 90] (-0.28,0.1);
 \end{tikzpicture}
}&=
\mathord{
\begin{tikzpicture}[baseline = -1mm]
     \node at (0.45,-0.25) {$\color{darkblue}\scriptstyle{r}$};
	\draw[-,thick,darkblue] (0.28,.4) to[out=240,in=90] (-0.28,-.3);
	\draw[<-,thick,darkblue] (-0.28,.4) to[out=300,in=90] (0.28,-0.3);
	\draw[-,thick,darkblue] (0.28,-0.3) to[out=-90, in=0] (0,-0.6);
	\draw[-,thick,darkblue] (0,-0.6) to[out = 180, in = -90] (-0.28,-0.3);
      \node at (0.27,-0.25) {$\dot$};
\end{tikzpicture}
}=0&\text{and}
\qquad\qquad
\mathord{
\begin{tikzpicture}[baseline = 2mm]
  \draw[->,thick,darkblue] (0.3,0.3) to[out=90,in=0] (0,.6);
  \draw[-,thick,darkblue] (0,0.6) to[out=180,in=90] (-.3,0.3);
\draw[-,thick,darkblue] (-.3,0.3) to[out=-90,in=180] (0,0);
  \draw[-,thick,darkblue] (0,0) to[out=0,in=-90] (0.3,0.3);
   \node at (0.3,0.3) {$\dot$};
   \node at (0.5,0.3) {$\color{darkblue}\scriptstyle{r}$};
\end{tikzpicture}
}&= 
\delta_{r,-k-1}
1_\unit\label{life2}
\end{align}
if $0 \leq r < -k$.

\begin{lemma}
The following relations hold:
\begin{align}
\mathord{
\begin{tikzpicture}[baseline = 1.5mm]
	\draw[->,thick,darkblue] (0.3,.6) to (-0.3,-.2);
	\draw[->,thick,darkblue] (0.3,-.2) to (-0.3,.6);
      \node at (-0.16,-0.02) {$\dot$};
\end{tikzpicture}
}
-
\mathord{
\begin{tikzpicture}[baseline = 1.5mm]
	\draw[->,thick,darkblue] (0.3,.6) to (-0.3,-.2);
	\draw[->,thick,darkblue] (0.3,-.2) to (-0.3,.6);
     \node at (0.16,0.42) {$\dot$};
 \end{tikzpicture}
}
&=
\mathord{
\begin{tikzpicture}[baseline = 1.5mm]
 	\draw[<-,thick,darkblue] (-0.33,-.2) to [out=90,in=180] (-0.1,0.1);
 	\draw[-,thick,darkblue] (-0.1,0.1) to [out=0,in=90] (0.13,-.2);
	\draw[<-,thick,darkblue] (-0.33,.6) to [out=-90,in=180] (-0.1,.3);
	\draw[-,thick,darkblue] (-0.1,.3) to [out=0,in=-90](0.13,.6);
\end{tikzpicture}
},
&
\mathord{
\begin{tikzpicture}[baseline = 1.5mm]
	\draw[->,thick,darkblue] (0.3,.6) to (-0.3,-.2);
	\draw[->,thick,darkblue] (0.3,-.2) to (-0.3,.6);
      \node at (0.16,-0.02) {$\dot$};
\end{tikzpicture}
}
-\mathord{
\begin{tikzpicture}[baseline = 1.5mm]
 	\draw[->,thick,darkblue] (0.3,.6) to (-0.3,-.2);
	\draw[->,thick,darkblue] (0.3,-.2) to (-0.3,.6);
      \node at (-0.16,0.42) {$\dot$};
\end{tikzpicture}
}&=
\mathord{
\begin{tikzpicture}[baseline = 1.5mm]
 	\draw[<-,thick,darkblue] (-0.33,-.2) to [out=90,in=180] (-0.1,0.1);
 	\draw[-,thick,darkblue] (-0.1,0.1) to [out=0,in=90] (0.13,-.2);
	\draw[<-,thick,darkblue] (-0.33,.6) to [out=-90,in=180] (-0.1,.3);
	\draw[-,thick,darkblue] (-0.1,.3) to [out=0,in=-90](0.13,.6);
\end{tikzpicture}
},
\label{leftcross}\\\notag\\
\mathord{
\begin{tikzpicture}[baseline = 2mm]
	\draw[-,thick,darkblue] (0.4,0) to[out=90, in=0] (0.1,0.5);
	\draw[->,thick,darkblue] (0.1,0.5) to[out = 180, in = 90] (-0.2,0);
      \node at (0.38,0.2) {$\dot$};
\end{tikzpicture}
}
&=
\mathord{
\begin{tikzpicture}[baseline = 2mm]
	\draw[-,thick,darkblue] (0.4,0) to[out=90, in=0] (0.1,0.5);
	\draw[->,thick,darkblue] (0.1,0.5) to[out = 180, in = 90] (-0.2,0);
      \node at (-0.17,0.2) {$\dot$};
\end{tikzpicture}
}\:,
&\mathord{
\begin{tikzpicture}[baseline = 2mm]
	\draw[-,thick,darkblue] (0.4,0.5) to[out=-90, in=0] (0.1,0);
	\draw[->,thick,darkblue] (0.1,0) to[out = 180, in = -90] (-0.2,0.5);
      \node at (-0.17,0.2) {$\dot$};
\end{tikzpicture}
}
&=
\mathord{
\begin{tikzpicture}[baseline = 2mm]
	\draw[-,thick,darkblue] (0.4,0.5) to[out=-90, in=0] (0.1,0);
	\draw[->,thick,darkblue] (0.1,0) to[out = 180, in = -90] (-0.2,0.5);
      \node at (0.38,0.2) {$\dot$};
\end{tikzpicture}
}.
\label{leftspade}
\end{align}
\end{lemma}

\proof
To prove (\ref{leftcross}),
take the first equation from (\ref{spit2}) describing how dots slide
past rightward crossings, vertically compose on top and bottom
with $t'$, then simplify using (\ref{c})--(\ref{life2}). 
For (\ref{leftspade}), it suffices to prove the first equation, 
since the latter then
follows on applying $\omega$ (recalling the heuristic for $\omega$
explained just before Lemma~\ref{omega}).
If $k < 0$
we vertically compose on the bottom with the
isomorphism
$\E\otimes \F \oplus 
\unit^{\oplus (-k)}
\stackrel{\sim}{\rightarrow}
\F\otimes \E$ from (\ref{invrel2}) to reduce to checking the following:
\begin{align*}
\mathord{

}
\text{ for all $0 \leq r < -k$.}
\end{align*}
To establish the first identity here, commute the dot past the
crossing on each side using (\ref{spit2}), then use the vanishing of the
curl from (\ref{life2}).
The second identity follows using (\ref{cough1}).
Finally, we must prove the first equation from (\ref{leftspade}) when $k
\geq 0$. In view of the definition of the leftward cap from (\ref{c}), we
must show equivalently that
\begin{align*}
\mathord{
%
}.
\end{align*}
To see this, use (\ref{leftcross}) to commute the bottom dot past the crossing, then appeal to (\ref{cough2}).
\endproof

We also
give meaning to {negatively dotted bubbles} by making the
following definitions for $r < 0$:
\begin{align}
\mathord{

}
&:=
\left\{
\begin{array}{cl}
\!\!\!\!\mathord{
\begin{tikzpicture}[baseline = 0mm]
  \draw[-,thick,darkblue] (0.3,0.1) to[out=90,in=0] (0,.4);
  \draw[-,thick,darkblue] (0,0.4) to[out=180,in=90] (-.3,0.1);
\draw[->,thick,darkblue] (-.3,0.1) to[out=-90,in=180] (0,-0.2);
  \draw[-,thick,darkblue] (0,-0.2) to[out=0,in=-90] (0.3,0.1);
   \node at (0.3,0.1) {$\dot$};
   \node at (0,0.4) {$\diamond$};
   \node at (0.62,0.1) {$\color{darkblue}\scriptstyle{-k}$};
   \node at (0,0.6) {$\color{darkblue}\scriptstyle{-r-1}$};
\end{tikzpicture}
}&\text{if $r > k-1$,}\\
-1_\unit&\text{if $r=k-1$},\\
0
&\text{if $r< k-1$},
\end{array}%
\right.
\label{ig}
\end{align}

\begin{lemma}\label{infgrass}
The infinite Grassmannian relations 
from Theorem~\ref{thm2}(i) all hold.
\end{lemma}

\proof
The equation (\ref{fed1})
is implied by (\ref{life1})--(\ref{life2}) and (\ref{ig}).
For (\ref{fed3}), we may assume using $\omega$
that $k\geq 0$. When $t=0$ the result 
follows trivially using (\ref{fed1}).
When $t > 0$ we
have:
\begin{align*}
\sum_{\substack{r,s \in \Z \\ r+s=t-2}}
%
}
\right) \stackrel{(\ref{posneg})}{=} 0.
\end{align*}
This implies (\ref{fed3}).
\endproof

The next lemma expresses the decorated leftward cups and caps in terms
of the undecorated ones. It means 
that we will not need to use the diamond
notation again after this.

\begin{lemma}\label{advances}
The following holds:
\begin{align}\label{invA}
\begin{array}{ll}
\mathord{
\begin{tikzpicture}[baseline = 1mm]
	\draw[-,thick,darkblue] (0.4,0.4) to[out=-90, in=0] (0.1,0);
	\draw[->,thick,darkblue] (0.1,0) to[out = 180, in = -90] (-0.2,0.4);
      \node at (0.12,-0.2) {$\color{darkblue}\scriptstyle{r}$};
      \node at (0.12,0.01) {$\diamond$};
\end{tikzpicture}
}
=
-\sum_{s \geq 0}
\mathord{
\begin{tikzpicture}[baseline = -1mm]
	\draw[-,thick,darkblue] (0.3,0.4) to[out=-90, in=0] (0,0);
	\draw[->,thick,darkblue] (0,0) to[out = 180, in = -90] (-0.3,0.4);
  \node at (-0.25,0.15) {$\dot$};
   \node at (-0.45,0.15) {$\color{darkblue}\scriptstyle{s}$};
\end{tikzpicture}
}
\mathord{
\begin{tikzpicture}[baseline = 2mm]
  \draw[->,thick,darkblue] (0.3,0.3) to[out=90,in=0] (0,.6);
  \draw[-,thick,darkblue] (0,0.6) to[out=180,in=90] (-.3,0.3);
\draw[-,thick,darkblue] (-.3,0.3) to[out=-90,in=180] (0,0);
  \draw[-,thick,darkblue] (0,0) to[out=0,in=-90] (0.3,0.3);
   \node at (0.3,0.3) {$\dot$};
   \node at (.9,0.3) {$\color{darkblue}\scriptstyle{-r-s-2}$};
\end{tikzpicture}
}
&\text{if $0 \leq r < k$,}\\
\end{array}
\end{align}
\end{lemma}

\proof
We explain the first equality; the second may then 
be deduced by applying $\omega$ using also
(\ref{leftspade}).
Remembering the definition (\ref{cutting1}), 
it suffices to show on replacing each
$\mathord{
%
}
$
that the matrix product
$$
\left[
\mathord{
%
}\right]
$$
is the $(k+1)\times(k+1)$ identity matrix.
This may be checked quite routinely using (\ref{posneg})--(\ref{life1}) and
Lemma~\ref{infgrass};
cf. the proof of \cite[Corollary 3.3]{B} for a similar argument.
\endproof

If we substitute the formulae from Lemma~\ref{advances} into (\ref{posneg}),
we obtain:
\begin{align}\label{pos2}
\mathord{
%
}.\label{neg2}
\end{align}

\begin{lemma}\label{obstacle}
The curl relations from Theorem~\ref{thm2}(iv) all hold.
\end{lemma}

\proof
In the next paragraph, we will establish the following:
\begin{align}\label{everything}
\mathord{
%
}.
\end{align}
Then to obtain the curl relations in the form (\ref{dog1}), take the
dotted curls on the left hand side of those relations, use (\ref{hecke}) 
to commute the dots past the
upward crossing,
convert the crossing to a rightward one
using (\ref{rightadj}) and the definition of $t$, then apply (\ref{everything}).

For (\ref{everything}),
we first prove the first equation when $k \geq 0$:
\begin{align*}
\mathord{
%
}.
\end{align*}
The first equation when $k < 0$ is immediate from (\ref{life2}).
Then the second equation then follows
by applying $\omega$ and using (\ref{cough1}).
\endproof

The proofs of the next two lemmas are intertwined with each other.

\begin{lemma}\label{pitch}
The following relations hold:
\begin{align}
\label{final1}
\mathord{
%
}.\label{final2}
\end{align}
\end{lemma}

\proof
It suffices to prove the left hand equalities in
(\ref{final1})--(\ref{final2});
then the right hand ones follow by applying $\omega$.
In the next two paragraphs, we will prove the left hand equality in
(\ref{final1}) assuming $k
\leq 0$ and the left hand equality in (\ref{final2}) assuming $k > 0$.

Consider (\ref{final1}) when $k \leq 0$.
We claim that
\begin{equation}
\mathord{
%
}\,.\label{quarry}
\end{equation}
To prove this, vertically compose on the bottom with the
isomorphism
$$
\left[\:\:
\mathord{
%
}
\right]:\E\otimes \E \otimes \F \oplus \E^{\oplus (-k)}
\stackrel{\sim}{\rightarrow} \E\otimes \F\otimes \E
$$
to reduce to showing equivalently that
\begin{align*}
\mathord{
%
}\quad\text{ for $0 \leq r < -k$}.
\end{align*}
Here are the proofs of these two identities:
\begin{align*}
\mathord{
%
}.
\end{align*}
Thus, the claim (\ref{quarry}) is proved.
Then we have that
\begin{align*}
\mathord{
%
},\end{align*}
establishing (\ref{final1}).

Next consider (\ref{final2}) when $k > 0$.
The strategy to prove this is the same as in the previous paragraph.
One first verifies that
$\mathord{
%
}
$
by vertically composing on the top with the isomorphism
$$
\left[\:
\mathord{
%
}\:\:\right]^T:\E\otimes \F\otimes \E\stackrel{\sim}{\rightarrow} \F\otimes
\E\otimes \E \oplus \E^{\oplus k}.
$$
Then this can be used to show
$\mathord{
%
}.$

The partial results established so far are all that are needed to
prove Lemma~\ref{adjdone} below.
To complete the proof of the present lemma,
suppose first that $k > 0$.
We take the left hand equality from
(\ref{final2}) proved in the previous paragraph,
attach
leftward caps to the top left and top right strands, then
simplify using the left adjunction relations
to be established in Lemma~\ref{adjdone}.
This establishes (\ref{final1}) for $k > 0$.
Finally, (\ref{final2}) for $k \leq 0$ may be deduced from
(\ref{final1}) by a similar procedure.
\endproof

\begin{lemma}\label{adjdone}
The left adjunction relations from Theorem~\ref{thm2}(ii) hold.
\end{lemma}

\proof
As usual, it suffices to prove the first equality.
If $k \leq 0$ then
\begin{align*}
\mathord{
\begin{tikzpicture}[baseline = 0]
  \draw[-,thick,darkblue] (0.3,0) to (0.3,-.4);
	\draw[-,thick,darkblue] (0.3,0) to[out=90, in=0] (0.1,0.4);
	\draw[-,thick,darkblue] (0.1,0.4) to[out = 180, in = 90] (-0.1,0);
	\draw[-,thick,darkblue] (-0.1,0) to[out=-90, in=0] (-0.3,-0.4);
	\draw[-,thick,darkblue] (-0.3,-0.4) to[out = 180, in =-90] (-0.5,0);
  \draw[->,thick,darkblue] (-0.5,0) to (-0.5,.4);
\end{tikzpicture}
}
&\stackrel{(\ref{c})}{=}
\mathord{
\begin{tikzpicture}[baseline = 0]
  \draw[->,thick,darkblue] (-0.3,0.6) to (-0.3,.8);
	\draw[-,thick,darkblue] (-0.3,0.6) to[out=-90, in=90] (0.1,-0.3);
	\draw[-,thick,darkblue] (0.1,-0.3) to[out = -90, in = 0] (-0.1,-0.6);
	\draw[-,thick,darkblue] (-0.1,-0.6) to[out = 180, in = -90] (-0.3,-0.3);
	\draw[-,thick,darkblue] (-0.3,-0.3) to[out = 90, in = 180] (0.2,0.4);
	\draw[-,thick,darkblue] (0.2,0.4) to[out = 0, in =90] (0.5,0);
  \draw[-,thick,darkblue] (0.5,0) to (0.5,-.8);
      \node at (0.28,-0.3) {$\color{darkblue}\scriptstyle{-k}$};
      \node at (0.1,-0.3) {$\dot$};
\end{tikzpicture}
}
\stackrel{(\ref{final1})}{=}
\mathord{
\begin{tikzpicture}[baseline = 0]
  \draw[->,thick,darkblue] (0.5,0.6) to (0.5,.8);
	\draw[-,thick,darkblue] (0.5,0.6) to[out=-90, in=90] (0.1,-0.3);
	\draw[-,thick,darkblue] (0.1,-0.3) to[out = -90, in = 0] (-0.1,-0.6);
	\draw[-,thick,darkblue] (-0.1,-0.6) to[out = 180, in = -90] (-0.3,-0.3);
	\draw[-,thick,darkblue] (-0.3,-0.3) to[out = 90, in = 180] (0.2,0.4);
	\draw[-,thick,darkblue] (0.2,0.4) to[out = 0, in =90] (0.5,0);
  \draw[-,thick,darkblue] (0.5,0) to (0.5,-.8);
      \node at (.28,-0.3) {$\color{darkblue}\scriptstyle{-k}$};
      \node at (0.1,-0.3) {$\dot$};
\end{tikzpicture}
}
\:\:\substack{(\ref{dog1})\\{\displaystyle=}\\(\ref{fed1})}\:\:
\mathord{\begin{tikzpicture}[baseline=0]
  \draw[->,thick,darkblue] (0,-0.4) to (0,.4);
\end{tikzpicture}
}\,.\end{align*}
If $k > 0$ then
\begin{align*}
\mathord{
\begin{tikzpicture}[baseline = 0]
  \draw[-,thick,darkblue] (0.3,0) to (0.3,-.4);
	\draw[-,thick,darkblue] (0.3,0) to[out=90, in=0] (0.1,0.4);
	\draw[-,thick,darkblue] (0.1,0.4) to[out = 180, in = 90] (-0.1,0);
	\draw[-,thick,darkblue] (-0.1,0) to[out=-90, in=0] (-0.3,-0.4);
	\draw[-,thick,darkblue] (-0.3,-0.4) to[out = 180, in =-90] (-0.5,0);
  \draw[->,thick,darkblue] (-0.5,0) to (-0.5,.4);
\end{tikzpicture}
}
&\stackrel{(\ref{c})}{=}
\mathord{
\begin{tikzpicture}[baseline = 0]
  \draw[-,thick,darkblue] (0.3,-0.6) to (0.3,-.8);
	\draw[-,thick,darkblue] (0.3,-0.6) to[out=90, in=-90] (-0.1,0.3);
	\draw[-,thick,darkblue] (-0.1,0.3) to[out = 90, in = 180] (0.1,0.6);
	\draw[-,thick,darkblue] (0.1,0.6) to[out = 0, in = 90] (0.3,0.3);
	\draw[-,thick,darkblue] (0.3,0.3) to[out = -90, in = 0] (-0.2,-0.4);
	\draw[-,thick,darkblue] (-0.2,-0.4) to[out = 180, in =-90] (-0.5,0);
  \draw[->,thick,darkblue] (-0.5,0) to (-0.5,.8);
      \node at (-0.26,0.2) {$\color{darkblue}\scriptstyle{k}$};
      \node at (-0.1,0.2) {$\dot$};
\end{tikzpicture}
}
\stackrel{(\ref{final2})}{=}
\mathord{
\begin{tikzpicture}[baseline = 0]
  \draw[-,thick,darkblue] (-0.5,-0.6) to (-0.5,-.8);
	\draw[-,thick,darkblue] (-0.5,-0.6) to[out=90, in=-90] (-0.1,0.3);
	\draw[-,thick,darkblue] (-0.1,0.3) to[out = 90, in = 180] (0.1,0.6);
	\draw[-,thick,darkblue] (0.1,0.6) to[out = 0, in = 90] (0.3,0.3);
	\draw[-,thick,darkblue] (0.3,0.3) to[out = -90, in = 0] (-0.2,-0.4);
	\draw[-,thick,darkblue] (-0.2,-0.4) to[out = 180, in =-90] (-0.5,0);
  \draw[->,thick,darkblue] (-0.5,0) to (-0.5,.8);
      \node at (-0.26,0.2) {$\color{darkblue}\scriptstyle{k}$};
      \node at (-0.1,0.2) {$\dot$};
\end{tikzpicture}
}
\:\:\substack{(\ref{dog1})\\{\displaystyle=}\\(\ref{fed1})}\:\:
\mathord{\begin{tikzpicture}[baseline=0]
  \draw[->,thick,darkblue] (0,-0.4) to (0,.4);
\end{tikzpicture}
}\:.
\end{align*}
Note we have only used the parts of Lemma~\ref{pitch} that were already
proved without forward reference to the present lemma.
\endproof

There are just two more relations to be checked; the arguments here
are analogous to ones in \cite[$\S$3.1.2]{KL3}.

\begin{lemma}\label{bubbly}
The bubble slide relations from
Theorem~\ref{thm2}(v)
hold.
\end{lemma}

\proof
We just explain the argument for $k \geq 0$; the case $k < 0$ is similar.
We first prove (\ref{bs2}). This is trivial for $r < 0$ due to (\ref{fed1}), so we may
assume
that $r \geq 0$.
Then we calculate:
\begin{align*}
\mathord{\begin{tikzpicture}[baseline=0]
  \draw[<-,thick,darkblue] (0.3,0) to[out=90,in=0] (0,0.3);
  \draw[-,thick,darkblue] (0,0.3) to[out=180,in=90] (-.3,0);
\draw[-,thick,darkblue] (-.3,0) to[out=-90,in=180] (0,-0.3);
  \draw[-,thick,darkblue] (0,-0.3) to[out=0,in=-90] (0.3,0);
  \draw[->,thick,darkblue] (0.6,-0.5) to (0.6,.6);
   \node at (-.5,0) {$\color{darkblue}\scriptstyle{r}$};
      \node at (-.3,0) {$\dot$};
\end{tikzpicture}
}
&\stackrel{(\ref{neg})}{=}
\mathord{\begin{tikzpicture}[baseline=0]
  \draw[-,thick,darkblue] (0.4,-0.5) to [out=135,in=-90](0.1,.05);
  \draw[->,thick,darkblue] (0.1,0.05) to [out=90,in=-135] (0.4,.6);
  \draw[<-,thick,darkblue] (0.3,0) to[out=90,in=0] (0,0.3);
  \draw[-,thick,darkblue] (0,0.3) to[out=180,in=90] (-.3,0);
\draw[-,thick,darkblue] (-.3,0) to[out=-90,in=180] (0,-0.3);
  \draw[-,thick,darkblue] (0,-0.3) to[out=0,in=-90] (0.3,0);
  \node at (-.5,0) {$\color{darkblue}\scriptstyle{r}$};
      \node at (-.3,0) {$\dot$};
\end{tikzpicture}
}
\:\:\substack{(\ref{cough2})\\{\displaystyle =}\\(\ref{final2})}\:
\mathord{\begin{tikzpicture}[baseline=0]
  \draw[-,thick,darkblue] (-0.4,-0.5) to [out=45,in=-90](-.1,.05);
  \draw[->,thick,darkblue] (-0.1,0.05) to [out=90,in=-45] (-0.4,.6);
  \draw[<-,thick,darkblue] (0.3,0) to[out=90,in=0] (0,0.3);
  \draw[-,thick,darkblue] (0,0.3) to[out=180,in=90] (-.3,0);
\draw[-,thick,darkblue] (-.3,0) to[out=-90,in=180] (0,-0.3);
  \draw[-,thick,darkblue] (0,-0.3) to[out=0,in=-90] (0.3,0);
  \node at (-.5,0) {$\color{darkblue}\scriptstyle{r}$};
      \node at (-.3,0) {$\dot$};
\end{tikzpicture}
}
\stackrel{(\ref{hecke})}{=}
\mathord{\begin{tikzpicture}[baseline=0]
  \draw[-,thick,darkblue] (-0.4,-0.5) to [out=45,in=-90](-.1,.05);
  \draw[->,thick,darkblue] (-0.1,0.05) to [out=90,in=-45] (-0.4,.6);
  \draw[<-,thick,darkblue] (0.3,0) to[out=90,in=0] (0,0.3);
  \draw[-,thick,darkblue] (0,0.3) to[out=180,in=90] (-.3,0);
\draw[-,thick,darkblue] (-.3,0) to[out=-90,in=180] (0,-0.3);
  \draw[-,thick,darkblue] (0,-0.3) to[out=0,in=-90] (0.3,0);
  \node at (.25,.36) {$\color{darkblue}\scriptstyle{r}$};
      \node at (.1,.26) {$\dot$};
\end{tikzpicture}
}
+\sum_{\substack{s,t \geq 0 \\ s+t=r-1}}
\mathord{
\begin{tikzpicture}[baseline = -0.5mm]
	\draw[<-,thick,darkblue] (0,0.6) to (0,0.3);
	\draw[-,thick,darkblue] (0,0.3) to [out=-90,in=180] (.3,-0.2);
	\draw[-,thick,darkblue] (0.3,-0.2) to [out=0,in=-90](.5,0);
	\draw[-,thick,darkblue] (0.5,0) to [out=90,in=0](.3,0.2);
	\draw[-,thick,darkblue] (0.3,.2) to [out=180,in=90](0,-0.3);
	\draw[-,thick,darkblue] (0,-0.3) to (0,-0.6);
  \node at (-.2,.3) {$\color{darkblue}\scriptstyle{s}$};
      \node at (0,.3) {$\dot$};
   \node at (0.3,0.03) {$\color{darkblue}\scriptstyle{t}$};
      \node at (0.2,.17) {$\dot$};
\end{tikzpicture}
}\\
&\:\:\substack{(\ref{hecke})\\{\displaystyle =}\\(\ref{dog1})}\:
\mathord{\begin{tikzpicture}[baseline=0]
  \draw[<-,thick,darkblue] (0.3,0) to[out=90,in=0] (0,0.3);
  \draw[-,thick,darkblue] (0,0.3) to[out=180,in=90] (-.3,0);
\draw[-,thick,darkblue] (-.3,0) to[out=-90,in=180] (0,-0.3);
  \draw[-,thick,darkblue] (0,-0.3) to[out=0,in=-90] (0.3,0);
  \draw[->,thick,darkblue] (-0.6,-0.5) to (-0.6,.6);
   \node at (-.1,0) {$\color{darkblue}\scriptstyle{r}$};
      \node at (-.3,0) {$\dot$};
\end{tikzpicture}
}\:
-\sum_{\substack{s,t\geq 0 \\ s+t=r-1}} \sum_{m \geq 0}
\mathord{\begin{tikzpicture}[baseline=0]
  \draw[<-,thick,darkblue] (1.2,0) to[out=90,in=0] (0.9,0.3);
  \draw[-,thick,darkblue] (0.9,0.3) to[out=180,in=90] (.6,0);
\draw[-,thick,darkblue] (.6,0) to[out=-90,in=180] (0.9,-0.3);
  \draw[-,thick,darkblue] (0.8,-0.3) to[out=0,in=-90] (1.2,0);
  \draw[->,thick,darkblue] (-0.6,-0.5) to (-0.6,.6);
   \node at (-1,0) {$\color{darkblue}\scriptstyle{s+m}$};
      \node at (-.6,0) {$\dot$};
   \node at (0.1,0) {$\color{darkblue}\scriptstyle{t-m-1}$};
      \node at (.6,0) {$\dot$};
\end{tikzpicture}
}\:.
\end{align*}
This easily simplifies to the right hand side of (\ref{bs2}).

Now we deduce (\ref{bs1}). Let $u$ be an indeterminant
and 
\begin{align}\label{symser}
e(u) &:= \sum_{r \geq 0} e_r u^{-r}, 
&h(u) &:= \sum_{r \geq 0} h_r
u^{-r}
\end{align} 
be the generating functions for the elementary and complete symmetric
functions. These are elements of $\Sym[[u^{-1}]]$
which satisfy $e(u) h(-u) = 1$.
Lemma~\ref{infgrass} implies that the homomorphism $\beta$ defined after (\ref{savvy}) satisfies
\begin{align}\label{eser}
\beta(e(u)) &=
-\sum_{r \geq 0} 
\mathord{
\begin{tikzpicture}[baseline = 1.25mm]
  \draw[<-,thick,darkblue] (0,0.4) to[out=180,in=90] (-.2,0.2);
  \draw[-,thick,darkblue] (0.2,0.2) to[out=90,in=0] (0,.4);
 \draw[-,thick,darkblue] (-.2,0.2) to[out=-90,in=180] (0,0);
  \draw[-,thick,darkblue] (0,0) to[out=0,in=-90] (0.2,0.2);
   \node at (-0.2,0.2) {$\dot$};
   \node at (-0.7,0.2) {$\color{darkblue}\scriptstyle{r+k-1}$};
\end{tikzpicture}}\:
u^{-r},\\\label{hser}
\beta(h(-u))& = \sum_{r \geq 0} 
\mathord{
\begin{tikzpicture}[baseline = 1.25mm]
  \draw[->,thick,darkblue] (0.2,0.2) to[out=90,in=0] (0,.4);
  \draw[-,thick,darkblue] (0,0.4) to[out=180,in=90] (-.2,0.2);
\draw[-,thick,darkblue] (-.2,0.2) to[out=-90,in=180] (0,0);
  \draw[-,thick,darkblue] (0,0) to[out=0,in=-90] (0.2,0.2);
   \node at (0.2,0.2) {$\dot$};
   \node at (0.7,0.2) {$\color{darkblue}\scriptstyle{r-k-1}$};
\end{tikzpicture}
}u^{-r}.
\end{align}
Also let
$p(u) := \sum_{r \geq 0} (r+1) x^r u^{-r-2}$,
where $x$ is the upward dot as usual.
The identity (\ref{bs2}) just proved asserts that $$
\beta(e(u))\otimes 
1_\E = 1_\E \otimes \beta(e(u)) - p(u)
\otimes \beta(e(u)).
$$
Multiplying on the left and right by $\beta(h(-u)) = \beta(e(u))^{-1}$, 
we deduce that
$$
1_\E\otimes  \beta(h(-u)) = \beta(h(-u)) \otimes 1_\E -
\beta(h(-u)) \otimes p(u).
$$ 
This is equivalent to (\ref{bs1}).
\endproof

\begin{lemma}\label{altbraidprop}
The alternating braid relation from
Theorem~\ref{thm2}(vi) holds.
\end{lemma}

\proof
Again, we just sketch the argument when $k \geq 0$, since $k < 0$ is
similar.
The idea is to attach crossings to the top left and bottom right 
pairs of strands of the second equality of (\ref{spit1}) to deduce that
$$
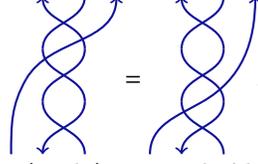
\mathord{
\begin{tikzpicture}[baseline = 0]
	\draw[->,thick,darkblue] (0.28,0.5) to[out=90, in=-90] (-.28,1.2);
	\draw[-,thick,darkblue] (0.28,1.2) to[out = -90, in = 90] (-0.28,0.5);
	\draw[-,thick,darkblue] (0.28,-0.2) to[out=90,in=-90] (-0.28,.5);
	\draw[-,thick,darkblue] (-0.28,-0.2) to[out=90,in=-90] (0.28,.5);
	\draw[-,thick,darkblue] (0.28,-.9) to[out=90,in=-90] (-0.28,-0.2);
	\draw[<-,thick,darkblue] (-0.28,-.9) to[out=90,in=-90] (0.28,-0.2);
	\draw[-,thick,darkblue] (-.7,-.9) to [out=90,in=-120] (-.5,0.1);
	\draw[->,thick,darkblue] (-.5,0.1) to [out=60,in=-90] (.7,1.2);
\end{tikzpicture}
}
=
\mathord{
\begin{tikzpicture}[baseline = 0]
	\draw[->,thick,darkblue] (0.28,0.5) to[out=90, in=-90] (-.28,1.2);
	\draw[-,thick,darkblue] (0.28,1.2) to[out = -90, in = 90] (-0.28,0.5);
	\draw[-,thick,darkblue] (0.28,-0.2) to[out=90,in=-90] (-0.28,.5);
	\draw[-,thick,darkblue] (-0.28,-0.2) to[out=90,in=-90] (0.28,.5);
	\draw[-,thick,darkblue] (0.28,-.9) to[out=90,in=-90] (-0.28,-0.2);
	\draw[<-,thick,darkblue] (-0.28,-.9) to[out=90,in=-90] (0.28,-0.2);
	\draw[-,thick,darkblue] (-.7,-.9) to [out=90,in=-150] (0,-.2);
	\draw[->,thick,darkblue] (0,-0.2) to [out=30,in=-90] (.7,1.2);
\end{tikzpicture}
}.
$$
Now apply (\ref{pos})--(\ref{neg}) to remove $t \circ t'$ and $t'
\circ t$ on each side then simplify; along the way many bubbles and
curls vanish thanks to (\ref{fed1}) and (\ref{life1}).
\endproof

\section{Proofs of Theorems}

\proof[Proof of Theorem~\ref{thm1}]
We first establish the existence of $c'$ and $d'$ satisfying the
relations (\ref{pos})--(\ref{leftcurl}). 
So let $\H_k$ be as in Definition~\ref{maindef}. Define $t'$ and the
decorated leftward cups and caps from (\ref{cutting1})--(\ref{cutting2}), then define
$c'$, $d'$ and the negatively dotted bubbles by
(\ref{c}) and (\ref{ig}). 
We need to show that this $t'$ and these negatively dotted bubbles are the
same as the ones defined in the statement of Theorem~\ref{thm1}.
For $t'$, this follows from (\ref{final1}) and the left adjunction
relations (\ref{adjfinal}) proved in Lemma~\ref{adjdone}.
For the negatively dotted bubbles, the infinite Grassmannian relations
(\ref{fed1})--(\ref{fed3}) 
proved in Lemma~\ref{infgrass} are all that are needed to construct the
homomorphism $\beta$ from (\ref{savvy}).
In the ring of symmetric functions, it is well known that
\begin{equation}\label{pieri}
h_r = \det\left(e_{i-j+1}\right)_{i,j=1,\dots,r}.
\end{equation}
Hence, applying the automorphism of $\Sym$ that interchanges $h_r$ and
$(-1)^r e_r$, we get also
$(-1)^re_r = \det\left((-1)^{i-j+1}h_{i-j+1}\right)_{i,j=1,\dots,r}$.
On applying $\beta$, this shows that
\begin{align*}
(-1)^r \mathord{
\begin{tikzpicture}[baseline = 1.25mm]
  \draw[->,thick,darkblue] (0.2,0.2) to[out=90,in=0] (0,.4);
  \draw[-,thick,darkblue] (0,0.4) to[out=180,in=90] (-.2,0.2);
\draw[-,thick,darkblue] (-.2,0.2) to[out=-90,in=180] (0,0);
  \draw[-,thick,darkblue] (0,0) to[out=0,in=-90] (0.2,0.2);
   \node at (0.2,0.2) {$\dot$};
   \node at (0.7,0.2) {$\color{darkblue}\scriptstyle{r-k-1}$};
\end{tikzpicture}}
&=
\det\left(-
\mathord{
\begin{tikzpicture}[baseline = 1.25mm]
  \draw[<-,thick,darkblue] (0,0.4) to[out=180,in=90] (-.2,0.2);
  \draw[-,thick,darkblue] (0.2,0.2) to[out=90,in=0] (0,.4);
 \draw[-,thick,darkblue] (-.2,0.2) to[out=-90,in=180] (0,0);
  \draw[-,thick,darkblue] (0,0) to[out=0,in=-90] (0.2,0.2);
   \node at (-0.2,0.2) {$\dot$};
   \node at (-0.65,0.2) {$\color{darkblue}\scriptstyle{i-j+k}$};
\end{tikzpicture}
}
\right)_{i,j=1,\dots,r},\\
-(-1)^r\mathord{
\begin{tikzpicture}[baseline = 1.25mm]
  \draw[<-,thick,darkblue] (0,0.4) to[out=180,in=90] (-.2,0.2);
  \draw[-,thick,darkblue] (0.2,0.2) to[out=90,in=0] (0,.4);
 \draw[-,thick,darkblue] (-.2,0.2) to[out=-90,in=180] (0,0);
  \draw[-,thick,darkblue] (0,0) to[out=0,in=-90] (0.2,0.2);
   \node at (-0.2,0.2) {$\dot$};
   \node at (-0.7,0.2) {$\color{darkblue}\scriptstyle{r+k-1}$};
\end{tikzpicture}
}&=
\det\left(\mathord{
\begin{tikzpicture}[baseline = 1.25mm]
  \draw[->,thick,darkblue] (0.2,0.2) to[out=90,in=0] (0,.4);
  \draw[-,thick,darkblue] (0,0.4) to[out=180,in=90] (-.2,.2);
\draw[-,thick,darkblue] (-.2,0.2) to[out=-90,in=180] (0,0);
  \draw[-,thick,darkblue] (0,0) to[out=0,in=-90] (0.2,0.2);
   \node at (0.2,0.2) {$\dot$};
   \node at (0.65,0.2) {$\color{darkblue}\scriptstyle{i-j-k}$};
\end{tikzpicture}
}\right)_{i,j=1,\dots,r},
\end{align*}
which easily simplify to produce the identities
(\ref{d1})--(\ref{d2}).
Thus, we are indeed in the setup of Theorem~\ref{thm1}. Now we get
the relations (\ref{pos})--(\ref{leftcurl}) from
(\ref{pos2})--(\ref{neg2}), the infinite Grassmannian relations
(\ref{fed1})--(\ref{fed3})
proved in Lemma~\ref{infgrass}, and
the curl relations (\ref{dog1}) proved in Lemma~\ref{obstacle}.

Next let $\mathcal{C}$ be a strict monoidal category with generators
$x,s,c,d,c',d'$ subject to the relations (\ref{hecke})--(\ref{rightadj}) and
(\ref{pos})--(\ref{leftcurl}). We have just demonstrated that all of these
relations hold in $\H_k$, hence, there is a strict $\k$-linear monoidal functor
$A:\mathcal{C} \rightarrow \H_k$ taking objects $\E,\F$ and 
generating morphisms $x,s,c,d,c',d'$ in $\mathcal{C}$ to the elements with the
same names in $\H_k$. 

In the other direction, we claim that there is a strict $\k$-linear
monoidal functor $B:\H_k\rightarrow \mathcal{C}$
sending the generating objects $\E, \F$ and morphisms
$x,s,c,d$ in $\H_k$ to the elements with the same
names in $\mathcal{C}$; this will eventually turn out to be a
two-sided inverse to $A$.
To prove the claim, we must verify that the three sets of defining relations of
$\H_k$ hold in
$\mathcal{C}$. It is immediate for (\ref{hecke}) and (\ref{rightadj}), so we are left with
checking the inversion relation.
We just do this in case $k \geq 0$, since the argument
for $k < 0$ is similar.
Defining the new morphisms
$$
\mathord{
\begin{tikzpicture}[baseline = 1mm]
	\draw[-,thick,darkblue] (0.4,0.4) to[out=-90, in=0] (0.1,0);
	\draw[->,thick,darkblue] (0.1,0) to[out = 180, in = -90] (-0.2,0.4);
      \node at (0.12,-0.2) {$\color{darkblue}\scriptstyle{r}$};
      \node at (0.12,0.01) {$\diamond$};
\end{tikzpicture}
}
:=
-\sum_{s \geq 0}
\mathord{
\begin{tikzpicture}[baseline = -1mm]
	\draw[-,thick,darkblue] (0.3,0.4) to[out=-90, in=0] (0,0);
	\draw[->,thick,darkblue] (0,0) to[out = 180, in = -90] (-0.3,0.4);
  \node at (-0.25,0.15) {$\dot$};
   \node at (-0.45,0.15) {$\color{darkblue}\scriptstyle{s}$};
\end{tikzpicture}
}
\mathord{
\begin{tikzpicture}[baseline = 2mm]
  \draw[->,thick,darkblue] (0.3,0.3) to[out=90,in=0] (0,.6);
  \draw[-,thick,darkblue] (0,0.6) to[out=180,in=90] (-.3,0.3);
\draw[-,thick,darkblue] (-.3,0.3) to[out=-90,in=180] (0,0);
  \draw[-,thick,darkblue] (0,0) to[out=0,in=-90] (0.3,0.3);
   \node at (0.3,0.3) {$\dot$};
   \node at (.86,0.3) {$\color{darkblue}\scriptstyle{-r-s-2}$};
\end{tikzpicture}
}
$$
in $\mathcal C$
for $r=0,1,\dots,k-1$, we claim that
$$
\left[
\displaystyle\mathord{
\begin{tikzpicture}[baseline = 0]
	\draw[->,thick,darkblue] (0.28,-.3) to (-0.28,.4);
	\draw[<-,thick,darkblue] (-0.28,-.3) to (0.28,.4);
\end{tikzpicture}
}\:\:\:
\mathord{
\begin{tikzpicture}[baseline = 1mm]
	\draw[-,thick,darkblue] (0.4,0.4) to[out=-90, in=0] (0.1,0);
	\draw[->,thick,darkblue] (0.1,0) to[out = 180, in = -90] (-0.2,0.4);
      \node at (0.12,-0.2) {$\color{darkblue}\scriptstyle{0}$};
      \node at (0.12,0.01) {$\diamond$};
\end{tikzpicture}
}
\:\:\:
\mathord{
\begin{tikzpicture}[baseline = 1mm]
	\draw[-,thick,darkblue] (0.4,0.4) to[out=-90, in=0] (0.1,0);
	\draw[->,thick,darkblue] (0.1,0) to[out = 180, in = -90] (-0.2,0.4);
      \node at (0.12,-0.2) {$\color{darkblue}\scriptstyle{1}$};
      \node at (0.12,0.01) {$\diamond$};
\end{tikzpicture}
}
\:\:\:\cdots\:\:\:
\mathord{
\begin{tikzpicture}[baseline = 1mm]
	\draw[-,thick,darkblue] (0.4,0.4) to[out=-90, in=0] (0.1,0);
	\draw[->,thick,darkblue] (0.1,0) to[out = 180, in = -90] (-0.2,0.4);
      \node at (0.12,-0.2) {$\color{darkblue}\scriptstyle{k-1}$};
      \node at (0.12,0.01) {$\diamond$};
\end{tikzpicture}
}
\right]
$$
is the the two-sided inverse of 
the morphism
(\ref{invrel1}).
Composing one way round gives the morphism
$$
\mathord{
\begin{tikzpicture}[baseline = 0mm]
	\draw[->,thick,darkblue] (0.28,0) to[out=90,in=-90] (-0.28,.6);
	\draw[-,thick,darkblue] (-0.28,0) to[out=90,in=-90] (0.28,.6);
	\draw[<-,thick,darkblue] (0.28,-.6) to[out=90,in=-90] (-0.28,0);
	\draw[-,thick,darkblue] (-0.28,-.6) to[out=90,in=-90] (0.28,0);
\end{tikzpicture}
}
-\sum_{r,s \geq 0}
\mathord{
\begin{tikzpicture}[baseline=-0.9mm]
	\draw[-,thick,darkblue] (0.3,0.6) to[out=-90, in=0] (0,0.1);
	\draw[->,thick,darkblue] (0,0.1) to[out = 180, in = -90] (-0.3,0.6);
      \node at (-0.4,0.3) {$\color{darkblue}\scriptstyle{r}$};
	\draw[<-,thick,darkblue] (0.3,-.6) to[out=90, in=0] (0,-0.1);
	\draw[-,thick,darkblue] (0,-0.1) to[out = 180, in = 90] (-0.3,-.6);
      \node at (-0.25,0.3) {$\dot$};
   \node at (-0.27,-0.4) {$\dot$};
   \node at (-.45,-.35) {$\color{darkblue}\scriptstyle{s}$};
\end{tikzpicture}}
\mathord{
\begin{tikzpicture}[baseline = 1mm]
  \draw[->,thick,darkblue] (0.2,0.2) to[out=90,in=0] (0,.4);
  \draw[-,thick,darkblue] (0,0.4) to[out=180,in=90] (-.2,0.2);
\draw[-,thick,darkblue] (-.2,0.2) to[out=-90,in=180] (0,0);
  \draw[-,thick,darkblue] (0,0) to[out=0,in=-90] (0.2,0.2);
   \node at (0.2,0.2) {$\dot$};
   \node at (.75,0.2) {$\color{darkblue}\scriptstyle{-r-s-2}$};
\end{tikzpicture}
},
$$
which is the identity by the relation (\ref{pos}) in $\mathcal C$.
The other way around, we get a $(k+1)\times (k+1)$-matrix. 
Its $1,1$-entry is the identity by (\ref{neg}). This is all that is
needed when $k=0$, but when $k > 0$ we also need 
to verifying the following for
$r,s=0,1,\dots,k-1$:
$$
\mathord{
\begin{tikzpicture}[baseline = -2mm]
	\draw[-,thick,darkblue] (0.28,-.6) to[out=120,in=-90] (-0.28,.1);
	\draw[<-,thick,darkblue] (-0.28,-.6) to[out=60,in=-90] (0.28,.1);
	\draw[-,thick,darkblue] (0.28,0.1) to[out=90, in=0] (0,0.4);
	\draw[-,thick,darkblue] (0,0.4) to[out = 180, in = 90] (-0.28,0.1);
      \node at (-0.45,0.05) {$\color{darkblue}\scriptstyle{r}$};
      \node at (-0.27,0.05) {$\dot$};
\end{tikzpicture}
}
=0,
\qquad
\mathord{
\begin{tikzpicture}[baseline = -2mm]
	\draw[<-,thick,darkblue] (0.28,.4) to[out=240,in=90] (-0.28,-.3);
	\draw[-,thick,darkblue] (-0.28,.4) to[out=300,in=90] (0.28,-0.3);
  	\draw[-,thick,darkblue] (0.28,-0.3) to[out=-90, in=0] (0,-0.6);
	\draw[-,thick,darkblue] (0,-0.6) to[out = 180, in = -90] (-0.28,-0.3);
      \node at (0,-0.8) {$\color{darkblue}\scriptstyle{s}$};
      \node at (0,-0.6) {$\diamond$};
\end{tikzpicture}
}=0,
\qquad
\mathord{
\begin{tikzpicture}[baseline = 2mm]
  \draw[<-,thick,darkblue] (0,0.6) to[out=180,in=90] (-.3,0.3);
  \draw[-,thick,darkblue] (0.3,0.3) to[out=90,in=0] (0,.6);
 \draw[-,thick,darkblue] (-.3,0.3) to[out=-90,in=180] (0,0);
  \draw[-,thick,darkblue] (0,0) to[out=0,in=-90] (0.3,0.3);
   \node at (-0.3,0.3) {$\dot$};
   \node at (-0.5,0.3) {$\color{darkblue}\scriptstyle{r}$};
      \node at (0,-0.2) {$\color{darkblue}\scriptstyle{s}$};
      \node at (0,0) {$\diamond$};
\end{tikzpicture}
}= \delta_{r,s} 1_\unit,
$$
Here is the proof of the first of these for $r=0,1,\dots,k-1$:
\begin{align}
\mathord{
\begin{tikzpicture}[baseline = -2mm]
	\draw[-,thick,darkblue] (0.28,-.6) to[out=120,in=-90] (-0.28,.1);
	\draw[<-,thick,darkblue] (-0.28,-.6) to[out=60,in=-90] (0.28,.1);
	\draw[-,thick,darkblue] (0.28,0.1) to[out=90, in=0] (0,0.4);
	\draw[-,thick,darkblue] (0,0.4) to[out = 180, in = 90] (-0.28,0.1);
      \node at (-0.45,0.05) {$\color{darkblue}\scriptstyle{r}$};
      \node at (-0.27,0.05) {$\dot$};
\end{tikzpicture}
}
&
\stackrel{(\ref{tprime})}{\displaystyle=}
\mathord{
\begin{tikzpicture}[baseline = -1mm]
	\draw[-,thick,darkblue] (0.25,.4) to [out=180,in=70] (-0.3,-.5);
        \draw[-,thick,darkblue] (0.5,0) to[out=90,in=0] (0.25,.4);
        \draw[-,thick,darkblue] (0.5,0) to[out=-90,in=0] (0.2,-.3);
        \draw[-,thick,darkblue] (-0.5,0.5) to[out=0,in=180] (0.2,-.3);
      \node at (0.13,0.21) {$\color{darkblue}\scriptstyle{r}$};
      \node at (-.05,0.21) {$\dot$};
        \draw[->,thick,darkblue] (-0.5,.5) to[out=180,in=-270] (-0.8,-.5);
\end{tikzpicture}
}
\stackrel{(\ref{hecke})}{\displaystyle=}
\mathord{
\begin{tikzpicture}[baseline = -1mm]
	\draw[-,thick,darkblue] (0.25,.4) to [out=180,in=70] (-0.3,-.5);
        \draw[-,thick,darkblue] (0.5,0) to[out=90,in=0] (0.25,.4);
        \draw[-,thick,darkblue] (0.5,0) to[out=-90,in=0] (0.2,-.3);
        \draw[-,thick,darkblue] (-0.5,0.5) to[out=0,in=180] (0.2,-.3);
      \node at (-0.4,-0.21) {$\color{darkblue}\scriptstyle{r}$};
      \node at (-.2,-0.21) {$\dot$};
        \draw[->,thick,darkblue] (-0.5,.5) to[out=180,in=-270] (-0.8,-.5);
\end{tikzpicture}
}-
\sum_{\substack{s,t\geq 0\\s+t=r-1}}
\mathord{
\begin{tikzpicture}[baseline = 0.5mm]
	\draw[-,thick,darkblue] (0.3,-.3) to[out=90, in=0] (0,0.6);
	\draw[->,thick,darkblue] (0,0.6) to[out = 180, in = 90] (-0.3,-0.3);
  \node at (0.3,0.15) {$\dot$};
   \node at (0.1,0.15) {$\color{darkblue}\scriptstyle{s}$};
\end{tikzpicture}
}\!\!
\mathord{
\begin{tikzpicture}[baseline = 2mm]
  \draw[-,thick,darkblue] (0.3,0.3) to[out=90,in=0] (0,.6);
  \draw[<-,thick,darkblue] (0,0.6) to[out=180,in=90] (-.3,0.3);
\draw[-,thick,darkblue] (-.3,0.3) to[out=-90,in=180] (0,0);
  \draw[-,thick,darkblue] (0,0) to[out=0,in=-90] (0.3,0.3);
   \node at (-0.3,0.3) {$\dot$};
   \node at (-.5,0.3) {$\color{darkblue}\scriptstyle{t}$};
\end{tikzpicture}
}
\stackrel{(\ref{rightcurl})}{\displaystyle=}0.\label{meeow}
\end{align}
To prove the second, note by definition for $s=0,1,\dots,k-1$ that
$$
\mathord{
\begin{tikzpicture}[baseline = -2mm]
	\draw[<-,thick,darkblue] (0.28,.4) to[out=240,in=90] (-0.28,-.3);
	\draw[-,thick,darkblue] (-0.28,.4) to[out=300,in=90] (0.28,-0.3);
  	\draw[-,thick,darkblue] (0.28,-0.3) to[out=-90, in=0] (0,-0.6);
	\draw[-,thick,darkblue] (0,-0.6) to[out = 180, in = -90] (-0.28,-0.3);
      \node at (0,-0.8) {$\color{darkblue}\scriptstyle{s}$};
      \node at (0,-0.6) {$\diamond$};
\end{tikzpicture}
}
=-
\sum_{r \geq 0}
\mathord{
\begin{tikzpicture}[baseline = -2mm]
	\draw[<-,thick,darkblue] (0.28,.4) to[out=240,in=90] (-0.28,-.3);
	\draw[-,thick,darkblue] (-0.28,.4) to[out=300,in=90] (0.28,-0.3);
  	\draw[-,thick,darkblue] (0.28,-0.3) to[out=-90, in=0] (0,-0.6);
	\draw[-,thick,darkblue] (0,-0.6) to[out = 180, in = -90] (-0.28,-0.3);
      \node at (-.48,-0.3) {$\color{darkblue}\scriptstyle{r}$};
      \node at (-.28,-0.3) {$\dot$};
\end{tikzpicture}
}\:\:
\mathord{
\begin{tikzpicture}[baseline = 2mm]
  \draw[->,thick,darkblue] (0.3,0.3) to[out=90,in=0] (0,.6);
  \draw[-,thick,darkblue] (0,0.6) to[out=180,in=90] (-.3,0.3);
\draw[-,thick,darkblue] (-.3,0.3) to[out=-90,in=180] (0,0);
  \draw[-,thick,darkblue] (0,0) to[out=0,in=-90] (0.3,0.3);
   \node at (0.3,0.3) {$\dot$};
   \node at (.9,0.3) {$\color{darkblue}\scriptstyle{-r-s-2}$};
\end{tikzpicture}
}.
$$
By the definition (\ref{d1}), the dotted bubble here is zero if $r \geq k$,
while for $r=0,1,\dots,k-1$ the dotted curl is zero by a similar
argument to (\ref{meeow}).
For the final relation involving the decorated dotted bubble, define
$\beta:\Sym \rightarrow \End_{\mathcal C}(\unit)$ by sending
$e_r \mapsto 
-\mathord{
\begin{tikzpicture}[baseline = 1.25mm]
  \draw[<-,thick,darkblue] (0,0.4) to[out=180,in=90] (-.2,0.2);
  \draw[-,thick,darkblue] (0.2,0.2) to[out=90,in=0] (0,.4);
 \draw[-,thick,darkblue] (-.2,0.2) to[out=-90,in=180] (0,0);
  \draw[-,thick,darkblue] (0,0) to[out=0,in=-90] (0.2,0.2);
   \node at (-0.2,0.2) {$\dot$};
   \node at (-0.7,0.2) {$\color{darkblue}\scriptstyle{r+k-1}$};
\end{tikzpicture}}$ for each $r \geq 0$.
Then by (\ref{pieri}), we have that $\beta((-1)^r h_r) = \det\left(\mathord{
\begin{tikzpicture}[baseline = 1.25mm]
  \draw[<-,thick,darkblue] (0,0.4) to[out=180,in=90] (-.2,0.2);
  \draw[-,thick,darkblue] (0.2,0.2) to[out=90,in=0] (0,.4);
 \draw[-,thick,darkblue] (-.2,0.2) to[out=-90,in=180] (0,0);
  \draw[-,thick,darkblue] (0,0) to[out=0,in=-90] (0.2,0.2);
   \node at (-0.2,0.2) {$\dot$};
   \node at (-0.7,0.2) {$\color{darkblue}\scriptstyle{i-j+k}$};
\end{tikzpicture}}\right)_{i,j=1,\dots,r}$. Assuming $r\leq k$, this is
exactly the definition 
of
$\mathord{
\begin{tikzpicture}[baseline = 1mm]
  \draw[->,thick,darkblue] (0.2,0.2) to[out=90,in=0] (0,.4);
  \draw[-,thick,darkblue] (0,0.4) to[out=180,in=90] (-.2,0.2);
\draw[-,thick,darkblue] (-.2,0.2) to[out=-90,in=180] (0,0);
  \draw[-,thick,darkblue] (0,0) to[out=0,in=-90] (0.2,0.2);
   \node at (0.2,0.2) {$\dot$};
   \node at (.7,0.2) {$\color{darkblue}\scriptstyle{r-k-1}$};
\end{tikzpicture}
}$
from 
(\ref{d1}).
Now suppose that $0 \leq r,s < k$. Applying $\beta$ to the symmetric
function identity $\sum_{t=0}^{k-s-1} (-1)^{k-s-1-t} e_{r-k+1+t} h_{k-s-1-t} =
\delta_{r,s}$ 
and using (\ref{rightcurl})
shows that
$-\sum_{t=0}^{k-s-1}
\mathord{
\begin{tikzpicture}[baseline = 1.25mm]
  \draw[<-,thick,darkblue] (0,0.4) to[out=180,in=90] (-.2,0.2);
  \draw[-,thick,darkblue] (0.2,0.2) to[out=90,in=0] (0,.4);
 \draw[-,thick,darkblue] (-.2,0.2) to[out=-90,in=180] (0,0);
  \draw[-,thick,darkblue] (0,0) to[out=0,in=-90] (0.2,0.2);
   \node at (-0.2,0.2) {$\dot$};
   \node at (-0.55,0.2) {$\color{darkblue}\scriptstyle{r+t}$};
\end{tikzpicture}}\:\:
\mathord{
\begin{tikzpicture}[baseline = 1mm]
  \draw[->,thick,darkblue] (0.2,0.2) to[out=90,in=0] (0,.4);
  \draw[-,thick,darkblue] (0,0.4) to[out=180,in=90] (-.2,0.2);
\draw[-,thick,darkblue] (-.2,0.2) to[out=-90,in=180] (0,0);
  \draw[-,thick,darkblue] (0,0) to[out=0,in=-90] (0.2,0.2);
   \node at (0.2,0.2) {$\dot$};
   \node at (.77,0.2) {$\color{darkblue}\scriptstyle{-s-t-2}$};
\end{tikzpicture}
}=\delta_{r,s}1_\unit$, which is exactly the identity we need.
This proves the claim,
so the functor $B$
is well-defined.

Next we check that $c'$ and $d'$ are the {\em unique} morphisms in $\mathcal
C$ satisfying the relations (\ref{pos})--(\ref{leftcurl}). We do this
by using the assumed relations to derive 
expressions for $c'$ and $d'$ in terms of the other generators. 
Note by the claim in the previous paragraph that the leftward
crossing $t'$ may be characterized as the first entry of the inverse
of the morphism (\ref{invrel1}) when $k \geq 0$; similarly, it is the
first entry of the inverse of the morphism (\ref{invrel2}) when $k <
0$. 
This shows that $t'$ does not depend
on the values of $c'$ and $d'$ (despite being defined in terms of them). 
Then, when $k \geq 0$, we argue as in (\ref{meeow})
to show that
$$
\mathord{
\begin{tikzpicture}[baseline = -2mm]
	\draw[-,thick,darkblue] (0.28,-.6) to[out=120,in=-90] (-0.28,.1);
	\draw[<-,thick,darkblue] (-0.28,-.6) to[out=60,in=-90] (0.28,.1);
	\draw[-,thick,darkblue] (0.28,0.1) to[out=90, in=0] (0,0.4);
	\draw[-,thick,darkblue] (0,0.4) to[out = 180, in = 90] (-0.28,0.1);
      \node at (-0.45,0.05) {$\color{darkblue}\scriptstyle{k}$};
      \node at (-0.27,0.05) {$\dot$};
\end{tikzpicture}
}
\stackrel{(\ref{tprime})}{\displaystyle=}
\mathord{
\begin{tikzpicture}[baseline = -1mm]
	\draw[-,thick,darkblue] (0.25,.4) to [out=180,in=70] (-0.3,-.5);
        \draw[-,thick,darkblue] (0.5,0) to[out=90,in=0] (0.25,.4);
        \draw[-,thick,darkblue] (0.5,0) to[out=-90,in=0] (0.2,-.3);
        \draw[-,thick,darkblue] (-0.5,0.5) to[out=0,in=180] (0.2,-.3);
      \node at (0.13,0.21) {$\color{darkblue}\scriptstyle{k}$};
      \node at (-.05,0.21) {$\dot$};
        \draw[->,thick,darkblue] (-0.5,.5) to[out=180,in=-270] (-0.8,-.5);
\end{tikzpicture}
}
\stackrel{(\ref{hecke})}{\displaystyle=}
\mathord{
\begin{tikzpicture}[baseline = -1mm]
	\draw[-,thick,darkblue] (0.25,.4) to [out=180,in=70] (-0.3,-.5);
        \draw[-,thick,darkblue] (0.5,0) to[out=90,in=0] (0.25,.4);
        \draw[-,thick,darkblue] (0.5,0) to[out=-90,in=0] (0.2,-.3);
        \draw[-,thick,darkblue] (-0.5,0.5) to[out=0,in=180] (0.2,-.3);
      \node at (-0.4,-0.21) {$\color{darkblue}\scriptstyle{k}$};
      \node at (-.2,-0.21) {$\dot$};
        \draw[->,thick,darkblue] (-0.5,.5) to[out=180,in=-270] (-0.8,-.5);
\end{tikzpicture}
}-
\sum_{\substack{s,t\geq 0\\s+t=k-1}}
\mathord{
\begin{tikzpicture}[baseline = 0.5mm]
	\draw[-,thick,darkblue] (0.3,-.3) to[out=90, in=0] (0,0.6);
	\draw[->,thick,darkblue] (0,0.6) to[out = 180, in = 90] (-0.3,-0.3);
  \node at (0.3,0.15) {$\dot$};
   \node at (0.1,0.15) {$\color{darkblue}\scriptstyle{s}$};
\end{tikzpicture}
}\!\!
\mathord{
\begin{tikzpicture}[baseline = 2mm]
  \draw[-,thick,darkblue] (0.3,0.3) to[out=90,in=0] (0,.6);
  \draw[<-,thick,darkblue] (0,0.6) to[out=180,in=90] (-.3,0.3);
\draw[-,thick,darkblue] (-.3,0.3) to[out=-90,in=180] (0,0);
  \draw[-,thick,darkblue] (0,0) to[out=0,in=-90] (0.3,0.3);
   \node at (-0.3,0.3) {$\dot$};
   \node at (-.5,0.3) {$\color{darkblue}\scriptstyle{t}$};
\end{tikzpicture}
}
\stackrel{(\ref{rightcurl})}{\displaystyle=}
\mathord{
\begin{tikzpicture}[baseline = 1mm]
	\draw[-,thick,darkblue] (0.4,0) to[out=90, in=0] (0.1,0.4);
	\draw[->,thick,darkblue] (0.1,0.4) to[out = 180, in = 90] (-0.2,0);
\end{tikzpicture}
}\:.
$$
This establishes the uniqueness of $d'$ when $k
\geq 0$.
Similarly, using (\ref{leftcurl}) in place of (\ref{rightcurl}), one
gets 
that
$$\mathord{
\begin{tikzpicture}[baseline = 1mm]
	\draw[-,thick,darkblue] (0.4,0.4) to[out=-90, in=0] (0.1,0);
	\draw[->,thick,darkblue] (0.1,0) to[out = 180, in = -90] (-0.2,0.4);
 \end{tikzpicture}
}\:=\:\mathord{
\begin{tikzpicture}[baseline = 0]
	\draw[-,thick,darkblue] (0.28,.6) to[out=240,in=90] (-0.28,-.1);
	\draw[<-,thick,darkblue] (-0.28,.6) to[out=300,in=90] (0.28,-0.1);
	\draw[-,thick,darkblue] (0.28,-0.1) to[out=-90, in=0] (0,-0.4);
	\draw[-,thick,darkblue] (0,-0.4) to[out = 180, in = -90] (-0.28,-0.1);
      \node at (0.54,-0.1) {$\color{darkblue}\scriptstyle{-k}$};
      \node at (0.27,-0.1) {$\dot$};
\end{tikzpicture}
}
$$
when $k \leq 0$, hence, 
$c'$ is unique when $k \leq 0$.
It remains to prove the uniqueness of $c'$ when $k > 0$ and of $d'$
when $k < 0$.
In the case that $k > 0$,
the claim from the previous paragraph shows
that 
the last entry of the inverse of (\ref{invrel1}) is
$$
-\sum_{s \geq 0}
\mathord{
\begin{tikzpicture}[baseline = -1mm]
	\draw[-,thick,darkblue] (0.3,0.4) to[out=-90, in=0] (0,0);
	\draw[->,thick,darkblue] (0,0) to[out = 180, in = -90] (-0.3,0.4);
  \node at (-0.25,0.15) {$\dot$};
   \node at (-0.45,0.15) {$\color{darkblue}\scriptstyle{s}$};
\end{tikzpicture}
}
\mathord{
\begin{tikzpicture}[baseline = 2mm]
  \draw[->,thick,darkblue] (0.3,0.3) to[out=90,in=0] (0,.6);
  \draw[-,thick,darkblue] (0,0.6) to[out=180,in=90] (-.3,0.3);
\draw[-,thick,darkblue] (-.3,0.3) to[out=-90,in=180] (0,0);
  \draw[-,thick,darkblue] (0,0) to[out=0,in=-90] (0.3,0.3);
   \node at (0.3,0.3) {$\dot$};
   \node at (.85,0.3) {$\color{darkblue}\scriptstyle{-s-k-1}$};
\end{tikzpicture}
}
\stackrel{(\ref{d1})}{\displaystyle=}
-\:
\mathord{
\begin{tikzpicture}[baseline = 1mm]
	\draw[-,thick,darkblue] (0.3,0.4) to[out=-90, in=0] (0,0);
	\draw[->,thick,darkblue] (0,0) to[out = 180, in = -90] (-0.3,0.4);
\end{tikzpicture}
}\,.
$$
Hence, $c'$ is unique when $k > 0$. The uniqueness of $d'$ when $k <
0$ is proved similarly.

Now we can complete the proof of the theorem.
First we show that $\mathcal C$ and $\H_k$ are isomorphic,
thereby establishing
 the equivalent presentation from the statement of the theorem.
To see this, we check that the functors $A$ and $B$ are two-sided inverses.
We have that $A \circ B = \id_{\H_k}$
obviously. To see that $B \circ A = \id_{\mathcal C}$,
it is clear that $B\circ A$ is the identity on the generating morphisms $x,s,c,d$, and
follows on the morphisms $c',d'$ by the uniqueness established in the
previous paragraph.
Finally, since $\H_k \cong \mathcal C$, the uniqueness of $c'$ and
$d'$ established in the previous paragraph implies they are also the
unique morphisms in $\H_k$ satisfying (\ref{pos})--(\ref{leftcurl}),
and we are done.
\endproof

\proof[Proof of Theorem~\ref{thm2}]
Parts (i), (ii), (iv), (v) and (vi) are proved in Lemmas~\ref{infgrass}, \ref{adjdone},
\ref{obstacle}, \ref{bubbly} and \ref{altbraidprop}, respectively.
Part (iii) for dots follows from (\ref{leftspade}), while for 
crossings it is an easy consequence of the ``pitchfork relations'' from
Lemma~\ref{pitch} (combined with the adjunction relations).
\endproof

\proof[Proof of Theorem~\ref{thm3}]
We first explain the identification with
Khovanov's category
 $\mathcal H$ from \cite[$\S$2.1]{K};
this
also follows from the more general 
identification with the Mackaay-Savage category made in the next
paragraph together with \cite[Remark 2.10]{MS} but it seems helpful
to treat this important special case independently.
So assume that $k=-1$.
Theorem~\ref{thm1} gives a presentation of $\H_{-1}$
with generating morphisms $x,s,c,d,c'$ and $d'$.
Comparing the relations
(\ref{hecke})--(\ref{rightadj}) and (\ref{pos})--(\ref{leftcurl}) with
the local relations in Khovanov's definition, we see that
there is a strict monoidal functor
$\H_{-1} \rightarrow \mathcal H$ sending $\E$ and $\F$ to Khovanov's objects
$\up=Q_+$ and $\down = Q_-$,
$s,c,d,c'$ and $d'$ to the
morphisms in Khovanov's category represented by the same diagrams, 
and $x$ to the right curl
$\mathord{
\begin{tikzpicture}[baseline = -1mm]
	\draw[<-,thick,darkblue] (0,0.4) to (0,0.15);
	\draw[-,thick,darkblue] (0,0.15) to [out=-90,in=180] (.15,-0.15);
	\draw[-,thick,darkblue] (0.15,-0.15) to [out=0,in=-90](.3,0);
	\draw[-,thick,darkblue] (0.3,0) to [out=90,in=0](.15,0.15);
	\draw[-,thick,darkblue] (0.15,0.15) to [out=180,in=90](0,-0.15);
	\draw[-,thick,darkblue] (0,-0.4) to (0,-0.15);
\end{tikzpicture}
}$.
This functor sends 
$\mathord{\begin{tikzpicture}[baseline = -1mm]
  \draw[-,thick,darkblue] (0,0.2) to[out=180,in=90] (-.2,0);
  \draw[-,thick,darkblue] (0.2,0) to[out=90,in=0] (0,.2);
 \draw[-,thick,darkblue] (-.2,0) to[out=-90,in=180] (0,-0.2);
  \draw[-,thick,darkblue] (0,-0.2) to[out=0,in=-90] (0.2,0);
 \end{tikzpicture}
}=\mathord{\begin{tikzpicture}[baseline = -1mm]
  \draw[-,thick,darkblue] (0,0.2) to[out=180,in=90] (-.2,0);
  \draw[->,thick,darkblue] (0.2,0) to[out=90,in=0] (0,.2);
 \draw[-,thick,darkblue] (-.2,0) to[out=-90,in=180] (0,-0.2);
  \draw[-,thick,darkblue] (0,-0.2) to[out=0,in=-90] (0.2,0);
      \node at (0.2,0) {$\dot$};
\end{tikzpicture}
}$
to the figure-of-eight, which is zero since it involves a left curl.
Hence, our functor factors through the specialization to induce
a functor from the additive envelope of $\H_{-1}(0)$ to $\mathcal H$.
To see that this functor is an isomorphism, we construct its two-sided
inverse.
This sends any diagram representing a morphism 
in Khovanov's category to the 
morphism in the additive envelope of 
$\H_{-1}(0)$ encoded by the same diagram.
It is well-defined since all of Khovanov's local relations hold in
$\H_{-1}(0)$, and also we have shown in Theorem~\ref{thm2} that $\H_{-1}(0)$ is strictly
pivotal
(something which
is required implicitly in Khovanov's definition).

For arbitrary $k \leq -1$, the identification of $\H_{k}(\delta)$ with the Mackaay-Savage
category
$\tilde{\mathcal H}^\lambda$ follows by a very similar argument.
Let $\lambda = \sum_{i} \lambda_i \omega_i$ be a dominant
weight (in the notation of \cite{MS}), 
and set  $k:=-\sum_i \lambda_i$
and $\delta := \sum_{i} i \lambda_i$.
In one direction,
the monoidal isomorphism from the additive envelope of
$\H_{k}(\delta)$ to $\tilde{\mathcal H}^\lambda$
sends 
our $x,s,c,d,c'$ and $d'$ to the morphisms in \cite{MS} denoted by the
same diagrams. The morphism denoted $c_n$ in \cite[(2.1)]{MS} 
for $0 \leq n \leq -k$ is our
$- \mathord{
\begin{tikzpicture}[baseline = 1.25mm]
  \draw[<-,thick,darkblue] (0,0.4) to[out=180,in=90] (-.2,0.2);
  \draw[-,thick,darkblue] (0.2,0.2) to[out=90,in=0] (0,.4);
 \draw[-,thick,darkblue] (-.2,0.2) to[out=-90,in=180] (0,0);
  \draw[-,thick,darkblue] (0,0) to[out=0,in=-90] (0.2,0.2);
   \node at (-0.2,0.2) {$\dot$};
   \node at (-0.75,0.2) {$\color{darkblue}\scriptstyle{n+k-1}$};
\end{tikzpicture}
}$, thanks to the definition of negatively dotted clockwise bubble at
the end of Theorem~\ref{thm1}.
Using this, it is straightforward to check that the local relations
in \cite[(2.2)--(2.9)]{MS} agree with the defining relations for
$\H_k(\delta)$ from 
(\ref{hecke})--(\ref{rightadj}) and (\ref{pos})--(\ref{leftcurl}).
Finally, $\H_k$ is strictly pivotal, which again is required
implicitly in the approach of \cite{MS}.
\endproof

\proof[Proof of Theorem~\ref{bt}]
By induction on the number of crossings, one checks using the relations established in $\S$2
that any diagram representing a morphism $\theta \in \Hom_{\H_k}(X,Y)$ can
be written as a $\Sym$-linear combination of morphisms in
$B_{\infty,\infty}(X,Y)$
with the same or fewer crossings.
So $B_{\infty,\infty}(X,Y)$ spans $\Hom_{\H_k}(X,Y)$. The problem is to prove it is
also linearly independent.
This is done already in the case $k=0$ in \cite[Theorem 1.2]{BCNR}. When $k <
0$, we will explain how to deduce it from \cite[Proposition 2.16]{MS}
in the next paragraph. Then it follows
for $k > 0$ by applying the isomorphism $\omega$ from
Lemma~\ref{omega}.

So assume henceforth that $k < 0$.
In order to make an observation about base change, let us add a superscript 
$\H_k^\k$ 
to indicate the ground ring: it suffices to establish linear independence for
$\H_k^\Z$; then one can obtain the linear independence
for arbitrary $\k$ by using the obvious functor
$\H^\k_k\rightarrow \H_k^\Z \otimes_{\Z}
\k$.
Thus we are reduced to the case that $\k = \Z$.
Suppose we are given some linear relation
$$
\sum_{\theta \in B_{\infty,\infty}(X,Y)} p_\theta \theta = 0
$$
for $p_\theta \in \Sym$.
Take any dominant integral weight $\lambda$ for $\mathfrak{sl}_\infty$
with $k=-\sum_i \lambda_i$, and set $\delta := \sum_i i \lambda_i$. By Theorem~\ref{thm3},
the specialized category $\H_k(\delta)$ embeds into the Mackaay-Savage
category $\tilde{\mathcal H}^\lambda$ 
over ground ring $\Z$. So we can appeal to \cite[Proposition 2.16]{MS}
to deduce
that $B_{\infty,\infty}(X,Y)$ is a basis for
$\Hom_{\H_k(\delta)}(X,Y)$
as a free right module over $\Sym$ specialized at $e_1 = -\delta$.
We deduce that $p_\theta|_{e_1 = -\delta} = 0$ for each $\theta$. Since there are infinitely many 
possibilities for $\delta$ as $\lambda$ varies (keeping $k < 0$ fixed), this is enough to show that all $p_\theta$ are
zero.
\endproof

\proof[Proof of Theorem~\ref{darn}]
Noting that $H_0^f \cong \k$, we denote the one-dimensional
$H_0^f$-module also by $\k$.
As $f(x_1) = 0$ in $H_1^f$, the functor $\ev \circ \Psi_f$ sends
$f(x)$ to zero, hence, it factors through the quotient category
$\H_{f,1}$ of $\H_k$.
Since $\k$ is a projective $H_0^f$-module and the induction and
restriction functors are biadjoint, it follows that $\ev \circ \Psi_f$
has image contained in the full subcategory $\bigoplus_{n \geq 0}
H_n^f \proj$ of $\bigoplus_{n \geq 0} H_n^f \Mod$.
This subcategory is additive and Karoubian, hence, 
the functor $\H_{f,1} \rightarrow \bigoplus_{n \geq 0} H_n^f \proj$
constructed so far extends to the functor $\psi_f$ on
$\Kar(\H_{f,1})$ from the statement of the theorem.

Now take $n \geq 0$. The functor $\psi_f$ 
maps $\E^{\otimes n}$ to $(\ind_{n-1}^n \circ \cdots \circ \ind_0^1) \k = H_n^f$, hence, it
defines an algebra homomorphism
\begin{equation}\label{tummy}
\psi_n:\End_{\H_{f,1}}(\E^{\otimes n})^\op \rightarrow
\End_{H_n^f}(H_n^f)^\op \equiv H_n^f.
\end{equation}
We claim that $\psi_n$ is actually an algebra isomorphism.
To see this, note by the relations that there is a homomorphism
\begin{align}\label{phin}
\phi_n:H^f_n
&\rightarrow \End_{\H_{f,1}}(\E^{\otimes n})^\op,\\\notag
x_i&\mapsto
(1_\E)^{\otimes(n-i)} \otimes x \otimes (1_\E)^{\otimes(i-1)},\\\notag
s_j&\mapsto (1_\E)^{\otimes(n-j-1)}\otimes s\otimes (1_\E)^{\otimes(j-1)}.
\end{align}
Now we observe that bubbles on the right edge are scalars in
$\End_{\H_{f,1}}(\E^{\otimes n})^\op$. This is straightforward to
prove directly at this point, but it also follows from the
more general statement made in the last part of
Lemma~\ref{belo}; the proof of that given below is independent of the
present theorem.
Hence, the easy spanning part of Theorem~\ref{bt}
implies that $\phi_n$ is surjective. 
Also $\psi_n \circ \phi_n = \id_{H_n^f}$ as the two sides agree on generators.
These two facts combined show that $\psi_n$ and $\phi_n$ are two-sided
inverses, and we have proved the claim.

By the claim, for any primitive idempotent
$e\in H_n^f$, there is a corresponding idempotent
$e\in \End_{\H_{f,1}}(\E^{\otimes n})$
defining an object $(\E^{\otimes n}, e) \in \Kar(\H_{f,1})$
which maps to $H_n^f e$ under the functor $\psi_f$.
This shows that the functor $\psi_f$ is dense.
It remains to show that it is full and faithful.
To see this, it suffices to take words $X = X_1\otimes\cdots\otimes
X_r$
and $Y=Y_1\otimes\cdots\otimes Y_r$ 
in the letters $\up$ and $\down$ representing objects of $\H_{f,1}$ 
such that
$$
n := \#\{i\:|\:X_i = \up\} - \#\{i\:|\:X_i = \down\}
=
\#\{j\:|\:Y_j = \up\} - \#\{j\:|\:Y_j = \down\},
$$
and show that 
$\psi_f:\Hom_{\H_{f,1}}(X,Y) \rightarrow \Hom_{H_n^f}(\psi_f(X),
\psi_f(Y))$
is an isomorphism.
To prove this, we first reduce to that case that $X = \unit$ using the
following commutative diagram, whose horizontal maps are the canonical isomorphisms
coming from adjunction/duality:
\begin{equation}\label{shrapnel}
\begin{CD}
\Hom_{\H_{f,1}}(X,Y) &@>\sim>>&\Hom_{\H_{f,1}}(\unit, X^* \otimes Y)\\
@V\psi_f VV&&@VV\psi_f V\\
\Hom_{H_n^f}(\psi_f(X), \psi_f(Y))
&@>\sim>>&\Hom_{H_0^f}(\k, \psi_f(X^* \otimes Y)).
\end{CD}
\end{equation}
Assume henceforth that $X = \unit$. We then proceed by induction on the length $s$
of $Y$, the case $s=0$ following since $\psi_0$ is an isomorphism.
If $s > 0$, then at least one letter $Y_i$ of $Y$ must equal $\down$.
If $i=s$, i.e., the letter $\down$ is on the right, then $Y \cong
\mathbf{0}$ as $1_\F = 0$ in $\H_{f,1}$, and the
conclusion is trivial. Otherwise, we may assume that $Y_i = \down$ and
$Y_{i+1}=\up$ for some $i < s$. Let $Y'$ be $Y$ with these two letters
interchanged and $Y''$ be $Y$ with these two letters removed.
Using the induction hypothesis and the following commutative diagram,
whose horizontal maps are the canonical isomorphisms coming from (\ref{invrel2}),
we see that the conclusion follows for $Y$ if we can prove it for
$Y'$:
\begin{equation}\label{bond}
\begin{CD}
\Hom_{\H_{f,1}}(\unit, Y)&@>\sim>>&\Hom_{\H_{f,1}}(\unit, Y' \oplus
Y''^{\oplus (-k)})\\
@V\psi_fVV&&@VV\psi_f V\\
\Hom_{H_0^f}(\k, \psi_f(Y))&@>\sim>>&\Hom_{H_0^f}(\unit, \psi_f(Y') \oplus \psi_f(Y'')^{\oplus{-k}}).
\end{CD}
\end{equation}
Repeating in this way, we can move the letter $\down$ of $Y$ to the right,
and then we are done as before.
\endproof

\proof[Proof of Lemma~\ref{belo}]
Suppose that \begin{align*}
f(u) &= u^\ell + z_1 u^{\ell-1}+\cdots + z_\ell,&
f'(u) &= u^{\ell'} + z_1' u^{\ell'-1}+\cdots + z_{\ell'}',
\end{align*}
for $z_1,\dots,z_\ell,z_1',\dots,z_{\ell'}' \in \k$.
Also set $z_0 = z_0' := 1$.

We first show that $\mathcal I_{f,f'}$ contains 
$\mathord{\begin{tikzpicture}[baseline = -1mm]
  \draw[-,thick,darkblue] (0,0.2) to[out=180,in=90] (-.2,0);
  \draw[->,thick,darkblue] (0.2,0) to[out=90,in=0] (0,.2);
 \draw[-,thick,darkblue] (-.2,0) to[out=-90,in=180] (0,-0.2);
  \draw[-,thick,darkblue] (0,-0.2) to[out=0,in=-90] (0.2,0);
   \node at (0.7,0) {$\color{darkblue}\scriptstyle{r-k-1}$};
      \node at (0.2,0) {$\dot$};
\end{tikzpicture}
}-\delta_r 1_\unit$
for all $r \geq 0$. Proceed by induction on $r$. If $r \leq \ell'$,
we are done by the definition of $\mathcal I_{f,f'}$, so assume that
$r > \ell'$. By (\ref{delta}), 
$u^k \delta(u) f(u) = f'(u)$, which is a polynomial in $u$. Hence, its
$u^{\ell'-r}$-coefficient is zero. This shows that
\begin{equation}\label{football}
\sum_{s=0}^{\ell} z_s \delta_{r-s} = 0.
\end{equation}
Since $r-k-1=\ell+r-\ell'-1 \geq \ell$, we can use 
$x^\ell + z_1 x^{\ell-1}+\cdots + z_\ell \in
\mathcal I_{f,f'}$ to deduce that
$
\sum_{s=0}^{\ell}
z_s \mathord{\begin{tikzpicture}[baseline = -1mm]
  \draw[-,thick,darkblue] (0,0.2) to[out=180,in=90] (-.2,0);
  \draw[->,thick,darkblue] (0.2,0) to[out=90,in=0] (0,.2);
 \draw[-,thick,darkblue] (-.2,0) to[out=-90,in=180] (0,-0.2);
  \draw[-,thick,darkblue] (0,-0.2) to[out=0,in=-90] (0.2,0);
   \node at (0.85,0) {$\color{darkblue}\scriptstyle{r-s-k-1}$};
      \node at (0.2,0) {$\dot$};
\end{tikzpicture}
}\in \mathcal I_{f,f'}.
$
Then by induction we get that
\begin{align*}
\mathord{\begin{tikzpicture}[baseline = -1mm]
  \draw[-,thick,darkblue] (0,0.2) to[out=180,in=90] (-.2,0);
  \draw[->,thick,darkblue] (0.2,0) to[out=90,in=0] (0,.2);
 \draw[-,thick,darkblue] (-.2,0) to[out=-90,in=180] (0,-0.2);
  \draw[-,thick,darkblue] (0,-0.2) to[out=0,in=-90] (0.2,0);
   \node at (0.7,0) {$\color{darkblue}\scriptstyle{r-k-1}$};
      \node at (0.2,0) {$\dot$};
\end{tikzpicture}
}
\!- \delta_r 1_\unit
 =
\mathord{\begin{tikzpicture}[baseline = -1mm]
  \draw[-,thick,darkblue] (0,0.2) to[out=180,in=90] (-.2,0);
  \draw[->,thick,darkblue] (0.2,0) to[out=90,in=0] (0,.2);
 \draw[-,thick,darkblue] (-.2,0) to[out=-90,in=180] (0,-0.2);
  \draw[-,thick,darkblue] (0,-0.2) to[out=0,in=-90] (0.2,0);
   \node at (0.7,0) {$\color{darkblue}\scriptstyle{r-k-1}$};
      \node at (0.2,0) {$\dot$};
\end{tikzpicture}
}
\!+ \sum_{s=1}^{\ell}
z_s \delta_{r-s} 1_\unit
\equiv
\sum_{s=0}^{\ell}
z_s \mathord{\begin{tikzpicture}[baseline = -1mm]
  \draw[-,thick,darkblue] (0,0.2) to[out=180,in=90] (-.2,0);
  \draw[->,thick,darkblue] (0.2,0) to[out=90,in=0] (0,.2);
 \draw[-,thick,darkblue] (-.2,0) to[out=-90,in=180] (0,-0.2);
  \draw[-,thick,darkblue] (0,-0.2) to[out=0,in=-90] (0.2,0);
   \node at (0.85,0) {$\color{darkblue}\scriptstyle{r-s-k-1}$};
      \node at (0.2,0) {$\dot$};
\end{tikzpicture}
}
\equiv 0
\end{align*}
modulo $\mathcal I_{f,f'}$, as required.

Next, let $e(u), h(u) \in \Sym[[u^{-1}]]$ be the power series from
(\ref{symser}).
The previous paragraph and (\ref{hser}) shows that $\beta(h(-u))
\equiv \delta(u) 1_\unit\pmod{\mathcal I_{f,f'}}$.
Since $e(u) = h(-u)^{-1}$ 
and $\delta'(u) = - \delta(u)^{-1}$, it follows that $\beta(e(u)) \equiv
-\delta'(u) 1_\unit\pmod{\mathcal I_{f,f'}}$.
In vew of (\ref{eser}), this shows that $\mathcal I_{f,f'}$ contains
$\mathord{\begin{tikzpicture}[baseline = -1mm]
  \draw[<-,thick,darkblue] (0,0.2) to[out=180,in=90] (-.2,0);
  \draw[-,thick,darkblue] (0.2,0) to[out=90,in=0] (0,.2);
 \draw[-,thick,darkblue] (-.2,0) to[out=-90,in=180] (0,-0.2);
  \draw[-,thick,darkblue] (0,-0.2) to[out=0,in=-90] (0.2,0);
   \node at (-0.72,0) {$\color{darkblue}\scriptstyle{r+k-1}$};
      \node at (-0.2,0) {$\dot$};
\end{tikzpicture}
}
-\delta_r' 1_\unit$ for all $r \geq 0$.

Now we can show that $f'(x') \in \mathcal I_{f,f'}$.
By (\ref{delta}), $z_r' = \sum_{s=0}^{\ell} z_s \delta_{r-s}$.
So
\begin{align*}
f'(x') &= \sum_{r=0}^{\ell'} z_r' 
\mathord{
\begin{tikzpicture}[baseline = 0]
	\draw[<-,thick,darkblue] (0.08,-.3) to (0.08,.4);
      \node at (0.08,0.07) {$\dot$};
      \node at (0.5,.07) {$\color{darkblue}\scriptstyle{\ell'-r}$};
\end{tikzpicture}
}
=
\sum_{r=0}^{\ell'} 
\sum_{s=0}^{\ell}
z_s
\delta_{r-s}
\mathord{
\begin{tikzpicture}[baseline = 0]
	\draw[<-,thick,darkblue] (0.08,-.3) to (0.08,.4);
      \node at (0.08,0.07) {$\dot$};
      \node at (0.5,.07) {$\color{darkblue}\scriptstyle{\ell'-r}$};
\end{tikzpicture}
}
\equiv
\sum_{s=0}^{\ell}
z_s
\sum_{r=0}^{\ell'} 
\mathord{
\begin{tikzpicture}[baseline = 0]
	\draw[<-,thick,darkblue] (0.08,-.3) to (0.08,.4);
      \node at (0.08,0.07) {$\dot$};
      \node at (0.5,.07) {$\color{darkblue}\scriptstyle{\ell'-r}$};
\end{tikzpicture}
}
\mathord{\begin{tikzpicture}[baseline = -1mm]
  \draw[-,thick,darkblue] (0,0.2) to[out=180,in=90] (-.2,0);
  \draw[->,thick,darkblue] (0.2,0) to[out=90,in=0] (0,.2);
 \draw[-,thick,darkblue] (-.2,0) to[out=-90,in=180] (0,-0.2);
  \draw[-,thick,darkblue] (0,-0.2) to[out=0,in=-90] (0.2,0);
   \node at (0.85,0) {$\color{darkblue}\scriptstyle{r-s-k-1}$};
      \node at (0.2,0) {$\dot$};
\end{tikzpicture}
}\\
&
\stackrel{(\ref{ig})}{\displaystyle=}
\sum_{s=0}^{\ell}
z_s
\sum_{r\geq 0}
\mathord{
\begin{tikzpicture}[baseline = 0]
	\draw[<-,thick,darkblue] (0.08,-.3) to (0.08,.4);
      \node at (0.08,0.07) {$\dot$};
      \node at (0.3,.07) {$\color{darkblue}\scriptstyle{r}$};
\end{tikzpicture}
}
\mathord{\begin{tikzpicture}[baseline = -1mm]
  \draw[-,thick,darkblue] (0,0.2) to[out=180,in=90] (-.2,0);
  \draw[->,thick,darkblue] (0.2,0) to[out=90,in=0] (0,.2);
 \draw[-,thick,darkblue] (-.2,0) to[out=-90,in=180] (0,-0.2);
  \draw[-,thick,darkblue] (0,-0.2) to[out=0,in=-90] (0.2,0);
   \node at (0.85,0) {$\color{darkblue}\scriptstyle{\ell-s-1-r}$};
      \node at (0.2,0) {$\dot$};
\end{tikzpicture}
}
\stackrel{(\ref{dog1})}{\displaystyle=}
\sum_{s=0}^{\ell}
z_s\:
\mathord{
\begin{tikzpicture}[baseline = -0.5mm]
	\draw[-,thick,darkblue] (0,0.6) to (0,0.3);
	\draw[-,thick,darkblue] (0,0.3) to [out=-90,in=180] (.3,-0.2);
	\draw[-,thick,darkblue] (0.3,-0.2) to [out=0,in=-90](.5,0);
	\draw[-,thick,darkblue] (0.5,0) to [out=90,in=0](.3,0.2);
	\draw[-,thick,darkblue] (0.3,.2) to [out=180,in=90](0,-0.3);
	\draw[->,thick,darkblue] (0,-0.3) to (0,-0.6);
   \node at (0.85,0.02) {$\color{darkblue}\scriptstyle{\ell-s}$};
      \node at (0.48,.02) {$\dot$};
\end{tikzpicture}
}\equiv 0\pmod{\mathcal I_{f,f'}}.
\end{align*}

So far, we have shown that the left tensor ideal generated by $f(x)$
and $\mathord{\begin{tikzpicture}[baseline = -1mm]
  \draw[-,thick,darkblue] (0,0.2) to[out=180,in=90] (-.2,0);
  \draw[->,thick,darkblue] (0.2,0) to[out=90,in=0] (0,.2);
 \draw[-,thick,darkblue] (-.2,0) to[out=-90,in=180] (0,-0.2);
  \draw[-,thick,darkblue] (0,-0.2) to[out=0,in=-90] (0.2,0);
   \node at (0.7,0) {$\color{darkblue}\scriptstyle{r-k-1}$};
      \node at (0.2,0) {$\dot$};
\end{tikzpicture}
}-\delta_r 1_\unit$ for $r=1,\dots,\ell'$
contains $f'(x')$ and $\mathord{\begin{tikzpicture}[baseline = -1mm]
  \draw[<-,thick,darkblue] (0,0.2) to[out=180,in=90] (-.2,0);
  \draw[-,thick,darkblue] (0.2,0) to[out=90,in=0] (0,.2);
 \draw[-,thick,darkblue] (-.2,0) to[out=-90,in=180] (0,-0.2);
  \draw[-,thick,darkblue] (0,-0.2) to[out=0,in=-90] (0.2,0);
   \node at (-0.7,0) {$\color{darkblue}\scriptstyle{r+k-1}$};
      \node at (-0.2,0) {$\dot$};
\end{tikzpicture}
}
-\delta_r' 1_\unit$ for $r=1,\dots\ell$.
Similar argument shows that the left tensor ideal generated by
the latter elements contains the former elements. This proves the lemma.
\endproof

\end{document}